\documentclass[12pt]{article}
\date{December 11, 2007 (revised February 29, 2008)}
\usepackage{amsthm,amsmath,amstext,amsfonts,amssymb}
\usepackage{color,graphics}
\usepackage[hypertex]{hyperref}
\newcommand{\Z}{{\mathbb Z}}
\newcommand{\R}{{\mathbb R}}
\renewcommand{\P}{{\mathbb P}}
\newcommand{\E}{{\mathbb E}}

\newcommand{\leb}{{\mathcal L}}
\newcommand{\ind}{{\mathbf 1}}
\newcommand{\red}{{\cal R}}
\newcommand{\blue}{{\cal B}}
\newcommand{\rred}{[{\cal R}]}
\newcommand{\bblue}{[{\cal B}]}
\newcommand{\mat}{{\cal M}}
\newcommand{\diam}{\mbox{\rm diam}}
\newcommand{\parend}{{\hfill $\Diamond$}\par\vskip2\parsep}
\newcommand{\seg}[1]{\langle #1 \rangle}
\newcommand{\sym}{\triangle}

\newtheorem{thm}{Theorem}
\newtheorem{lemma}[thm]{Lemma}
\newtheorem{prop}[thm]{Proposition}

\newcounter{mycount}
\newenvironment{mylist}{\begin{list}{{\rm (\roman{mycount})}}%
{\usecounter{mycount}\itemsep -3pt}}{\end{list}}

\newcommand{\dof}{\bf}

\title{Poisson Matching}

\author{Alexander E. Holroyd
\and Robin Pemantle \and Yuval Peres \and Oded Schramm }

\begin{document}
\maketitle \enlargethispage*{1.5cm}
\renewcommand{\thefootnote}{}
\footnotetext{{\bf\hspace{-6mm}Key words:} Poisson process, point
process, matching, stable marriage.}
\footnotetext{{\bf\hspace{-6mm}2000 Mathematics Subject
Classifications:} 60D05, 60G55, 05C70.}
\footnotetext{\hspace{-6mm}Funded by Microsoft Research (AEH, YP
\& OS), an NSERC discovery grant (AEH), NSF grants DMS-0603821
(RP) and DMS-0605166 (YP).}
\renewcommand{\thefootnote}{\arabic{footnote}}

\begin{abstract}
Suppose that red and blue points occur as independent homogeneous
Poisson processes in $\R^d$.  We investigate translation-invariant
schemes for perfectly matching the red points to the blue points.
For any such scheme in dimensions $d=1,2$, the matching distance
$X$ from a typical point to its partner must have infinite
$d/2$-th moment, while in dimensions $d\geq 3$ there exist schemes
where $X$ has finite exponential moments.  The Gale-Shapley stable
marriage is one natural matching scheme, obtained by iteratively
matching mutually closest pairs.  A principal result of this paper
is a power law upper bound on the matching distance $X$ for this
scheme. A power law lower bound holds also.  In particular, stable
marriage is close to optimal (in tail behavior) in $d=1$, but far
from optimal in $d\geq 3$.  For the problem of matching Poisson
points of a single color to each other, in $d=1$ there exist
schemes where $X$ has finite exponential moments, but if we insist
that the matching is a deterministic factor of the point process
then $X$ must have infinite mean.
%
%\bigskip
%
%Supposons que des points rouges et bleus \'evoluent suivant des processus de Poisson homog\`enes ind\'ependants dans $\R^d$. Nous nous int\'eressons \`a des proc\'ed\'es invariants par translation appariant de mani\`ere bijective les points rouges et les points bleus. En dimensions $d=1,2$, quelque soit le proc\'ed\'e consid\'er\'e, la distance d'appariement (matching distance) $X$ entre un point typique et son partenaire poss\`ede n\'ecessairement un $d/2$-\`eme moment infini. En revanche, en dimensions $d\geq 3$ il existe des proc\'ed\'es pour lesquels $X$ a des moments exponentiels finis. Le ``mariage stable'' de Gale-Shapley est un proc\'ed\'e naturel, obtenu en appariant une \`a une les paires mutuellement les plus proches. L'un des principaux r\'esultats de cet article est que dans le cas de ce proc\'ed\'e, la distance d'appariement $X$ est major\'ee par une loi de puissance. Une minoration en loi de puissance est \'egalement v\'erifi\'ee . En particulier, le mariage stable est essentiellement optimal (en terme de queue de distribution) en dimension $d=1$, mais il est loin d'\^etre optimal en dimensions $d\geq 3$. Dans le cas du probl\`eme  qui consiste \`a apparier des points d'une seule couleur issus d'un processus de Poisson, en dimension $d=1$ il existe des proc\'ed\'es pour lesquels $X$ a des moments exponentiels finis. Par contre, si l'on demande en plus que l'appariement soit une fonction d\'eterministe du processus ponctuel, alors $X$ a n\'ecessairement une moyenne infinie.
\end{abstract}

\section{Introduction}
%%%%%%%%%%%%%%%%%%%%%%

%
\begin{figure}
\addtolength{\hfuzz}{10cm}
 \noindent\hspace{-1.5cm}
\resizebox{8cm}{!}{\includegraphics{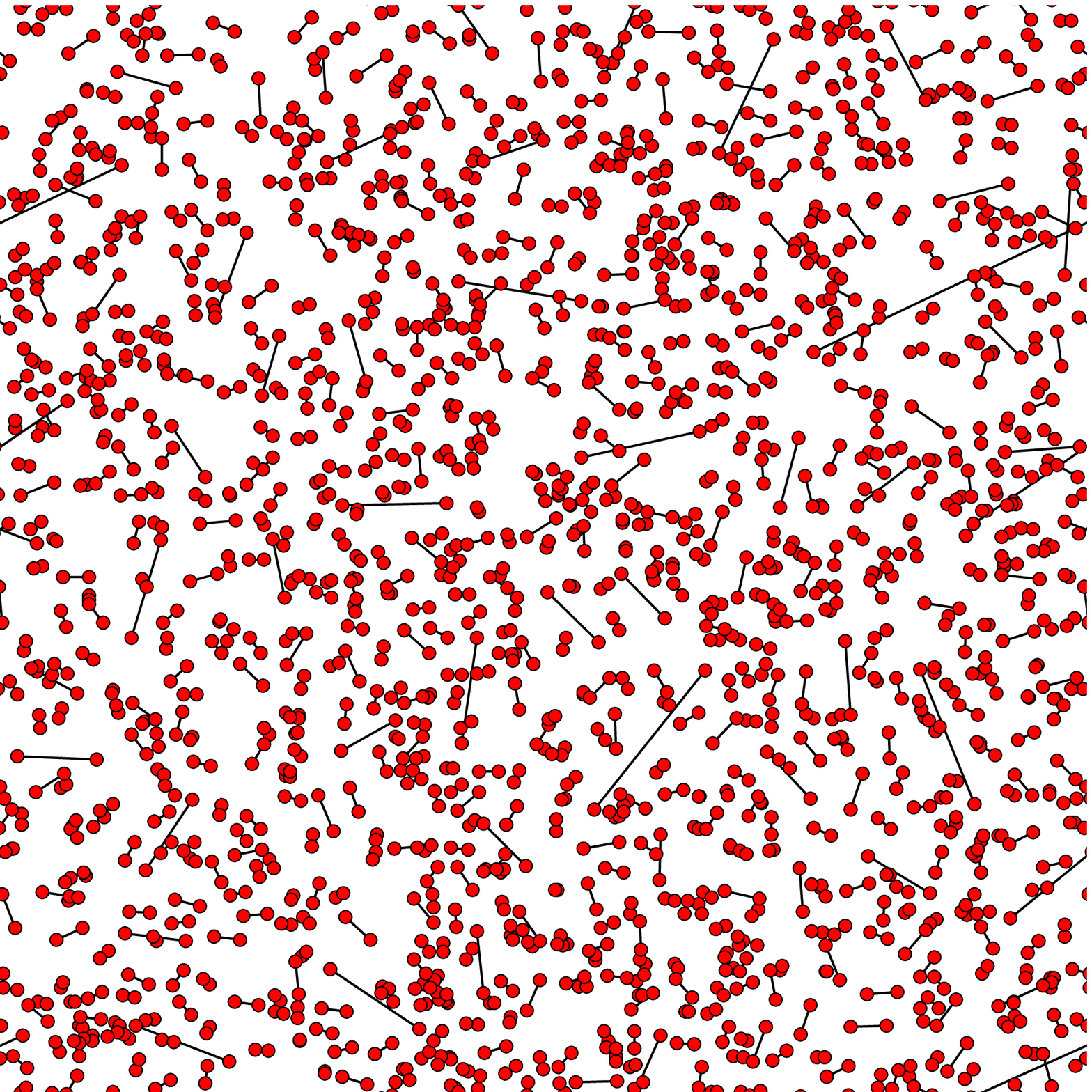}}
\resizebox{8cm}{!}{\includegraphics{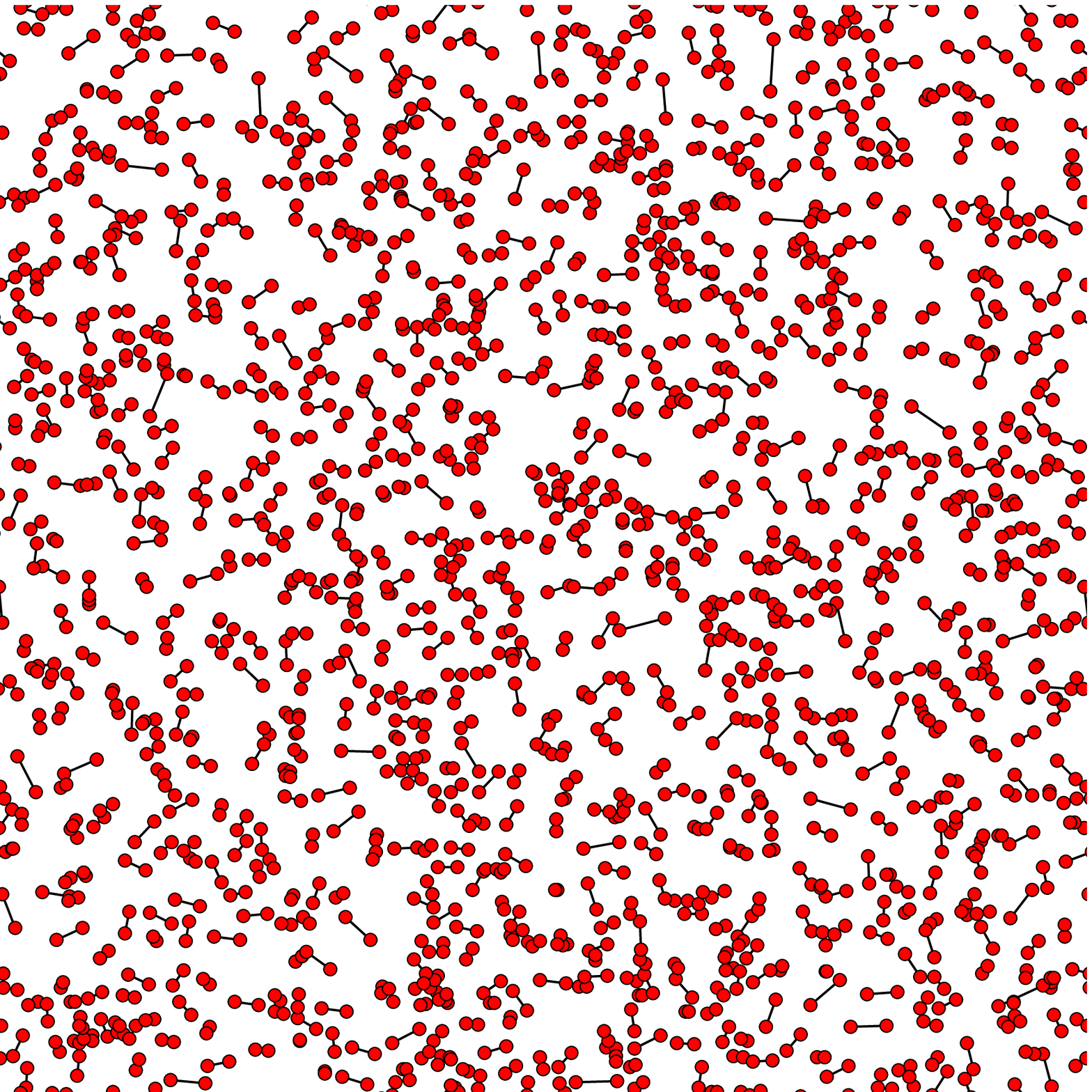}}
\\ \noindent\hspace*{-1.5cm}
\resizebox{8cm}{!}{\includegraphics{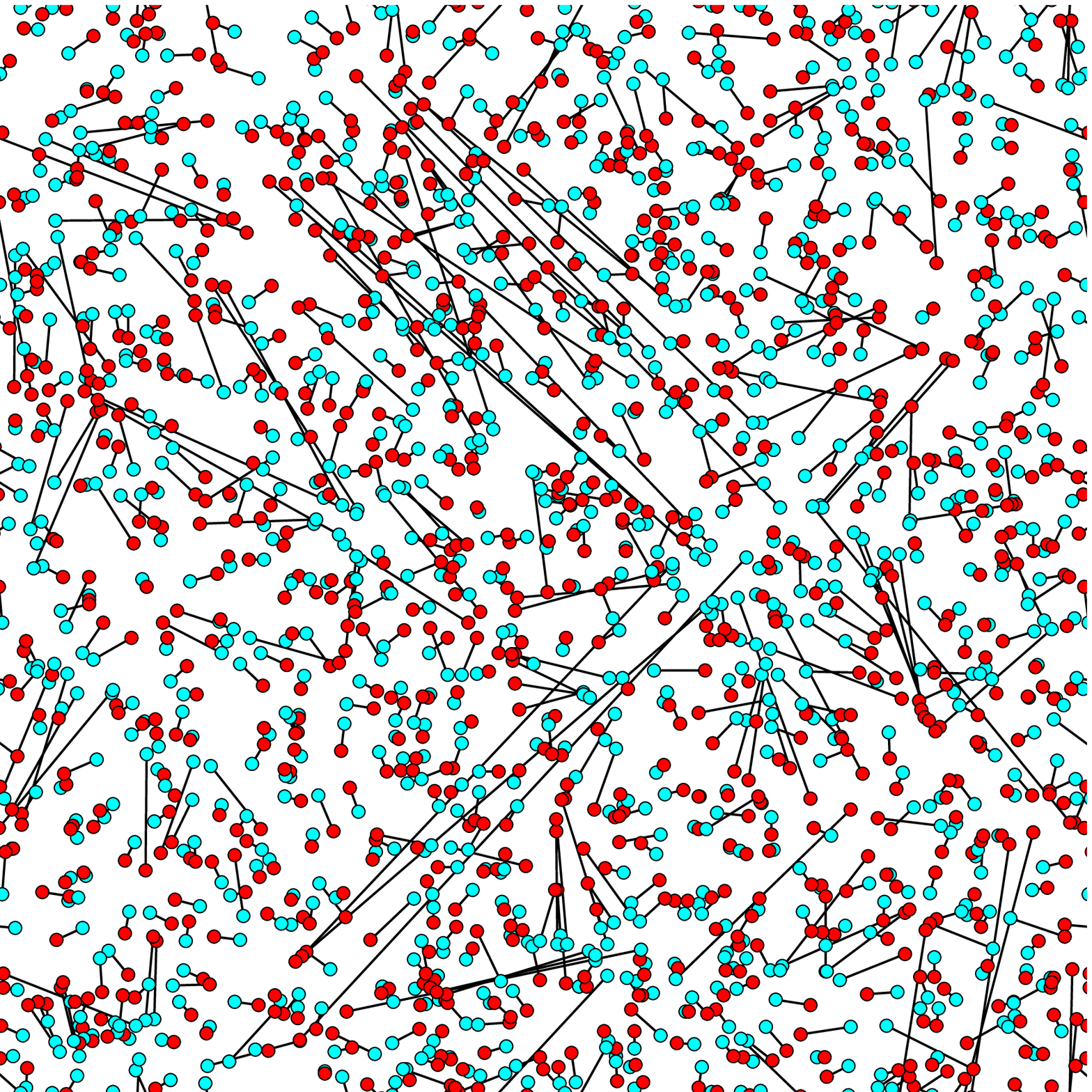}}
\resizebox{8cm}{!}{\includegraphics{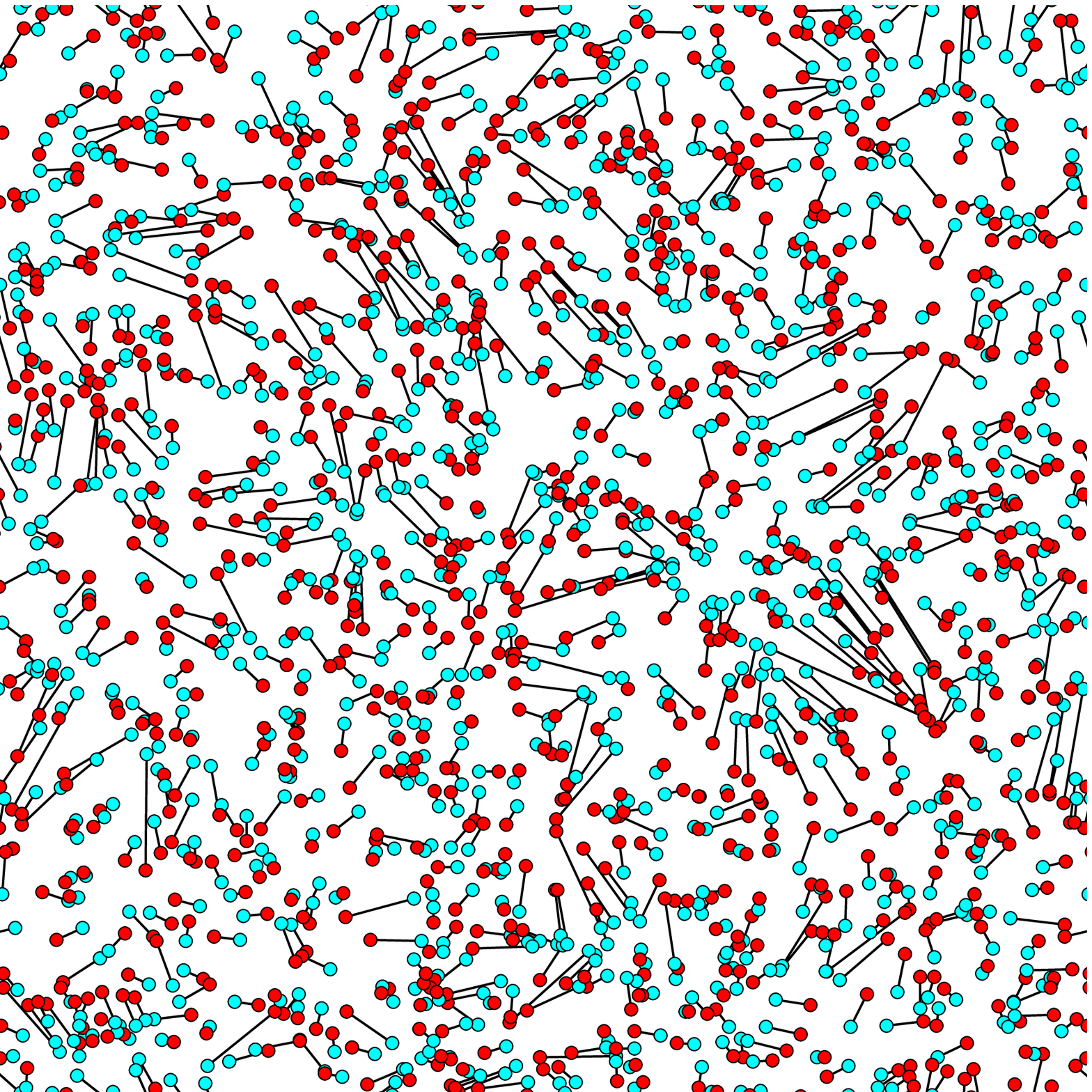}}

 \caption{Matchings of 2000 uniformly random points on a 2-dimensional torus:
(i) stable 1-color; (ii) minimum-length 1-color; (iii) stable
2-color; (iv) minimum-length 2-color.} \label{pics}
\addtolength{\hfuzz}{-10cm} \end{figure}

Let $\red$ be a simple point process of finite intensity in
$\R^d$. The {\dof support} of $\red$ is the random set
$\rred:=\{x\in \R^d: \red(\{x\})=1\}$.  Elements of $\rred$ are
called {\dof red points}. A {\dof one-color matching scheme} of
$\red$ is a simple point process $\mat$ of unordered pairs
$\{x,y\}\subset \R^d$, on a shared probability space, such that
almost surely $(V,E)=(\rred,[\mat])$ is a random graph which is a
perfect matching of $\rred$ (i.e.\ a simple graph with all degrees
1).  Let $\blue$ be a second simple point process, and call
elements of $\bblue$ {\dof blue points}. A {\dof two-color
matching scheme} between $\red$ and $\blue$ is a process $\mat$
which similarly yields almost surely a perfect bipartite matching
between $\rred$ and $\bblue$ (i.e.\ a perfect matching of
$\rred\cup\bblue$ where all the edges are from $\rred$ to
$\bblue$). In either case we denote by $\mat(x)$ the {\dof
partner} of a red or blue point $x$; that is the unique point such
that $\{x,\mat(x)\}\in[\mat]$.  See Figure \ref{pics} for some
examples on the finite torus.

We say that a one-color (respectively two-color) matching scheme
$\mat$ is {\dof translation-invariant} if the law of the joint process
$(\red,\mat)$ (respectively $(\red,\blue,\mat)$) is invariant
under translations of $\R^d$.   {\dof Isometry-invariance} is
defined analogously.  If almost surely $\mat=f(\red)$ (respectively
$\mat=f(\red,\blue)$) for some deterministic function $f$ then we
call $\mat$ a {\dof factor} matching scheme.  We sometimes refer to
a matching scheme which is not a factor as {\dof randomized}.

For a translation-invariant one-color or two-color matching scheme,
let $\P$ be the probability measure governing $(\red,\mat)$ or
$(\red,\blue,\mat)$, and $\E$ the associated expectation operator.
We are interested in the typical distance between matched pairs.
Assume without loss of generality that $\red$ has intensity $1$
(otherwise rescale). For $r\in [0,\infty]$ it is natural to consider
the quantity
$$F(r):=\E\#\big\{x\in \rred\cap [0,1)^d: |x-\mat(x)|\leq r\big\},$$
where $|\cdot|$ denotes the Euclidean norm. It is easy to see that
$F$ is a distribution function, therefore we introduce a random
variable $X$ with probability measure $\P^*$ and expectation
operator $\E^*$ such that
\begin{equation}
\label{def-x}
\P^*(X\leq r)=F(r).
\end{equation}
We can think of $X$ as the typical distance from a red point to its
partner.  In fact, this interpretation can be made rigorous via the
technology of Palm processes -- see Section \ref{prelim} below.

We consider the following main questions.  For Poisson processes on
$\R^d$, what is the best possible tail behavior (as measured by $X$)
for a translation-invariant (or isometry-invariant) matching scheme,
in the one-color and two-color cases? How do the answers depend on
dimension?  And if we insist on a factor matching scheme? We also
address the case of {\em stable} matchings -- see below.

Note the following trivial lower bound on $X$.  In the one-color or
two-color case, the partner of a point must be at least as far as
the closest other point.  In the case when $\red$ (respectively
$\red+\blue$) is a homogeneous Poisson process, this gives
\begin{equation}
\label{triv} \E^* e^{c X^d}=\infty
% \P^* (X>r) \geq e^{-c r^d} \quad\forall r>0,
\end{equation}
for some $c=c(d)\in(0,\infty)$.

The following theorems show that the optimal tail behavior of
two-color matching schemes depends dramatically on the dimension.
\begin{thm}[2-color upper bounds]
\label{2col-ub}
 \sloppy Let $\red,\blue$ be independent Poisson pro\-cess\-es of
intensity 1.  There exist isometry-invariant two-color matching
schemes satisfying: \fussy
\begin{mylist}
\item in $d=1$:\quad
$\P^*(X>r)\leq C r^{-1/2} \quad\forall r>0;$
\item in $d=2$:\quad
$\P^*(X>r)\leq C r^{-1} \quad\forall r>0;$
\item in $d\geq 3$:\quad
$\E^* e^{C X^d}<\infty.$
\end{mylist}
Here $C=C(d)\in(0,\infty)$ denotes a constant.  Furthermore, in (i)
the matching scheme is a factor.
\end{thm}

\begin{thm}[2-color lower bounds; \cite{h-p-extra}]
\label{2col-lb}
 Let $\red,\blue$ be independent Poisson processes of
intensity 1. In $d=1$ or $d=2$, any translation-invariant two-color
matching scheme (factor or not) satisfies
$$\E^* X^{d/2}=\infty.$$
\end{thm}

Together with the trivial bound \eqref{triv}, Theorems
\ref{2col-ub} and \ref{2col-lb} settle reasonably accurately the
question of optimal tail behavior for randomized two-color
matchings.  We do not know the optimal tail behavior for
translation-invariant {\em factor} matching schemes of two
independent Poisson processes in dimensions $d\geq 2$.
Theorem \ref{2col-lb} was
derived in \cite{h-p-extra} via results from
\cite{liggett-tagged,h-liggett}, the proofs of which were quite
involved.  We will present here a simple direct proof.

Here is a brief heuristic explanation for the above tail behavior,
and in particular the sharp difference between dimensions $d\leq
2$ and $d \geq 3$.  In a ball of large radius $r$, the discrepancy
between the numbers of red and blue points is typically a random
multiple of $r^{d/2}$ (by the central limit theorem).  This
discrepancy must be accommodated via the boundary of the ball,
which has size of order $r^{d-1}$.  When $d\leq 2$, the
discrepancy exceeds the boundary (with substantial probability),
so we expect that a fraction $r^{d/2}/r^d=r^{-d/2}$ of points must
be matched at distance at least of order $r$.  When $d\geq 3$, the
boundary exceeds the discrepancy, so far more efficient matching
is possible, with the tails determined by local deviations in the
distribution of points.

The following illustrates a case where allowing additional
randomization makes a striking difference to tail behavior.

\begin{thm}[1-color, 1 dimension]
\label{1col-1d} Let $d=1$, let $\red$ be a Poisson process of
intensity 1, and consider one-color matching schemes.
\begin{mylist}
\item Any translation-invariant {\em factor} matching satisfies
$\E^* X=\infty.$
\item There exists an isometry-invariant {\em randomized}
scheme satisfying $\P^*(X>r)=e^{-r} \quad\forall r>0.$
\end{mylist}
\end{thm}

The following shows that the above dichotomy does not extend to
higher dimensions.
\begin{thm}[1-color, 2 or more dimensions]
\label{trees}
 Let $\red$ be a Poisson process of intensity 1. For all
$d\geq 2$ there exists a translation-invariant one-color factor
matching scheme satisfying
$$\E^* e^{C X^d} <\infty,$$
for some $C=C(d)\in(0,\infty)$.  The same bound can be attained by a
randomized isometry-invariant matching scheme in all $d\geq 2$, and
by an isometry-invariant factor matching scheme in $d=2$.
\end{thm}
We do not know how to construct an isometry-invariant factor
matching satisfying the bound in Theorem \ref{trees} for $d\geq 3$
(see the remarks in Section \ref{sec-1-col}, however).

\paragraph{Stable matching: Iterated mutually closest matching algorithm.}
The following natural \lq\lq greedy\rq\rq\ algorithm gives a
matching scheme by trying to optimize locally. When considering a
two-color matching, call a pair of points $x,y$ {\dof potential
partners} if one is red while the other is blue. In the one-color
case, call $x$ and $y$ potential partners if they are distinct
points in $\rred$. We say that potential partners $x$ and $y$ are
{\dof mutually closest} if $y$ is the closest potential partner to
$x$ and $x$ is the closest potential partner to $y$. Now, given
the point configuration, match all mutually closest pairs to each
other, then remove these points and match all mutually closest
pairs in the remaining set of points. Repeat indefinitely.\parend

It turns out that the above algorithm yields a perfect matching
under general conditions (Proposition \ref{stable}), and in
particular this holds almost surely in the case when $\red$ is a
Poisson process (and for the two-color case, when $\blue$ is an
independent Poisson process of the same intensity). Furthermore, it
is the unique {\dof stable matching} in the sense of Gale and
Shapley \cite{gale-shapley}. (See Section \ref{prelim} for the
details). Evidently (under the aforementioned conditions) the stable
matching gives an isometry-invariant factor matching scheme.

We can accurately describe the tail behavior of stable matchings in
the one-color case, but some questions remain in the two-color case.

\begin{thm}[1-color stable matching]\label{stab-1col}
Let $\red$ be a Poisson process of intensity $1$.  For any  $d\geq
1$, the one-color stable matching satisfies
\begin{mylist}
\item\quad
$\E^* X^d=\infty;$
\item\quad
$\P^*(X>r)\leq C r^{-d} \quad{\forall r>0};$
\end{mylist}
for some $C=C(d)\in(0,\infty)$.
\end{thm}

\begin{thm}[2-color stable matching]\label{stab-2col}
Let $\red,\blue$ be independent Poisson processes of intensity $1$.
For $d\geq 1$, the two-color stable matching satisfies:
\begin{mylist}
\item\qquad $\E^* X^d=\infty$, \\(Theorem \ref{2col-lb} gives a better
bound in $d=1,2$);
 \item\qquad $\P^* (X>r)\leq C r^{-s} \quad{\forall r>0},$ \\
 where $s=s(d)\in(0,1)$ satisfies
 $s(1)=1/2$, and $C=C(d)\in(0,\infty)$.
\end{mylist}
\end{thm}
The power $s(d)$ is given explicitly as the solution of an equation,
and for example $s(2)=0.496\cdots$ and $s(3)=0.449\cdots$ -- see
Theorem \ref{t.upperbd}.
It is a fascinating unsolved question to determine the correct
power law for $d\geq 2$ (see the open problems at the end of the
article).

It is interesting that stable matching performs essentially
optimally (in terms of tail behavior) among (possibly randomized)
matching schemes in the two-color case for $d=1$, but not for
$d\geq 3$, and not in the one-color case.

\subsubsection*{Summary}

The following tables summarize the best known results for
isometry-invariant matchings of Poisson processes.

\vspace{5mm} \noindent
\begin{tabular}{|l|c|l|c|l|}
  \hline
\multicolumn{2}{|l|}{\bf 1-color matching} & Lower bound & & Upper
bound / \\
\multicolumn{2}{|l|}{} & &
 &  best construction
 \\
\hline\hline
 Randomized & All $d$ & $\E^*e^{cX^d}=\infty$ & $\approx$ & $\E^* e^{CX^{d}}<\infty$ \\
\hline
 Factor & $d=1$ & $\E^* X=\infty$ & $\approx$ & [$\P^*(X>r)\leq
 Cr^{-1}$]  \\
 & $d\geq 2$ & $\E^*e^{cX^d}=\infty$ & $\approx$ & $\E^* e^{CX^{d}}<\infty$ ((T) for $d\geq 3$)\\
\hline
 Stable & All $d$ & $\E^* X^d=\infty$ & $\approx$ & $\P^*(X>r)\leq Cr^{-d}$ \\
\hline
\end{tabular}

\vspace{5mm}

\noindent
\begin{tabular}{|l|c|l|c|l|}
  \hline
\multicolumn{2}{|l|}{\bf 2-color matching} &
 Lower bound & & Upper bound / \\
\multicolumn{2}{|l|}{} &
 & & best construction \\
\hline\hline
 Randomized & $d=1,2$ & $\E^* X^{d/2}=\infty$ & $\approx$ & $\P^*(X>r)\leq C r^{-d/2}$ \\
            & $d\geq 3$ & $\E^*e^{cX^d}=\infty$ & $\approx$ & $\E^* e^{CX^{d}}<\infty$ \\

\hline
 Factor & $d=1$ & [$\E^* X^{1/2}=\infty$] & $\approx$ & [$\P^*(X>r)\leq
 Cr^{-1/2}$] \\
                & $d=2$ & [$\E^* X=\infty$] & $\ll$ & $[\P^*(X>r)\leq Cr^{-0.496\cdots}]$ \\
                & $d\geq 3$ & $\E^*e^{cX^d}=\infty$ & $\ll$  & $[\P^*(X>r)\leq Cr^{-s(d)}]$ \\
\hline
Stable & $d=1$ & [$\E^* X^{1/2}=\infty$] & $\approx$ & $\P^*(X>r)\leq Cr^{-1/2}$ \\
        & $d=2$ & [$\E^* X=\infty$] & $\ll$ & $\P^*(X>r)\leq Cr^{-0.496\cdots}$ \\
        & $d\geq 3$ & $\E^* X^d=\infty$ & $\ll$ & $\P^*(X>r)\leq Cr^{-s(d)}$ \\
\hline
\end{tabular}

\vspace{5mm} \noindent
\begin{tabular}{|l|l|}
  \hline
{\bf Notes:} &\\
\hline
(T)& Translation-invariant scheme only. \\

[$\cdots$]& Bound follows from line above (lower bounds) \\ & \quad or below (upper bounds). \\

$\approx$& Indicates reasonably close lower and upper bounds.\\

$\ll$& Indicates a substantial gap between the lower \\ & \quad and upper bounds. \\

\hline
\end{tabular}

\subsubsection*{Extensions to other processes}

The results on 2-color matchings in Theorems \ref{2col-ub},
\ref{2col-lb}, \ref{stab-2col} all extend, with similar proofs, to
the following two variant settings.

\begin{mylist}
\item Perfect matchings of Heads to Tails for i.i.d.\ fair coin flips
indexed by $\Z^d$ (see \cite{h-p-extra,soo,timar-matching} for details).  (In order
to define stable matchings in this context, one must specify a way
to break ties between pairs of sites $\Z^d$ which are the same
distance apart).  Here the random variable $X$ denotes the distance
from the origin to its partner.

\item Fair allocations of Lebesgue measure to a Poisson process (see
\cite{hhp} and also \cite{grav,hhp-2,h-p-extra} for definitions and
background). Here $X$ denotes the distance from the origin to its
allocated Poisson point. The resulting upper bounds on the stable
allocation in Theorem \ref{stab-2col}(ii) represent a considerable
improvement on the previous best results: $X$ was known to have
finite $(1/18)$th moment in $d=1$, and no quantitative upper bound
was known in $d\geq 2$ \cite{hhp-2}.
\end{mylist}

Where appropriate, we make remarks following the proofs
regarding the adaptation of our results to these settings.
\medskip

Some of our results extend easily to point processes other than the
Poisson process. In particular, Theorem \ref{stab-1col}(ii) holds
for any translation-invariant simple point process for which the
stable matching is well-defined (see Proposition \ref{stable}).
Theorems \ref{stab-1col}(i) and \ref{stab-2col}(i) hold provided in
addition the point processes are tolerant of local modifications
(see Proposition \ref{local-mod}).

Another interesting variant concerns matching of random points on
large finite boxes; see e.g.\
\cite{akt,avis-davis-steele,frieze-mcdiarmid-reed,talagrand}.  We do
not explore this connection in depth, but our proof of Theorem
\ref{2col-ub}(iii) relies on the remarkable results of
\cite{talagrand}.

\section{Preliminaries}
%%%%%%%%%%%%%%%%%%%%%%%
\label{prelim}

In this section we present some useful elementary definitions and
results.

\subsubsection*{Some notation}

Let $\leb$ denote Lebesgue measure on $\R^d$, and denote the
Euclidean ball $B(x,r)$ $:=\{z\in\R^d:|x-z|<r\}$.  We denote the
unit cube $Q:=[0,1)^d\subset\R^d$, and $Q_u:=Q+u$ for $u\in\Z^d$.

\subsubsection*{Palm processes}

Consider a translation-invariant one-color or two-color matching
scheme, and let $\P$ be the probability measure governing
$(\red,\mat)$ or $(\red,\blue,\mat)$. We introduce the {\dof Palm
process} $(\red^*,\mat^*)$ or $(\red^*,\blue^*,\mat^*)$, with law
$\P^*$ and expectation $\E^*$, in which we condition on the presence
of a red point at the origin, while taking $\mat$ and (in the two
color-case) $\blue$ as a stationary background.  See e.g.\
\cite[Ch.\ 11]{kallenberg} for details.  In the case when $\red$ is
a homogeneous Poisson process, it turns out that $\red^*$ has the
same distribution as $\red$ with an {\em added} point at the origin:
\begin{equation}\label{poisson-palm}
[\red^*]\stackrel{d}{=}\rred\cup \{0\}.
\end{equation}
Also, if $\red$ and $\blue$ are independent processes then $\red^*$
and $\blue^*$ are independent, and $\blue^*\stackrel{d}{=}\blue$.

If $\red$ is a translation-invariant measure-valued process of
intensity $\lambda\in(0,\infty)$, and $(\red^*,\Psi^*)$ is the Palm
version of $\red$ taken together with any jointly
translation-invariant random background $\Psi$ (which may be a random
function or a random measure on $\R^d$), then the following properties
are standard (see \cite[Ch.\ 11]{kallenberg}).  Let $\theta^{x}$
denote translation by $x\in\R^d$ (defined to act on measures via
$(\theta^x\pi)(S)=\pi(S-x)$ for $S\subseteq\R^d$).  For any measurable
$S\subseteq\R^d$ and any event $A$ we have
\begin{equation}\label{palm1}
\E\int_S \ind\big[\theta^{-x}(\red,\Psi)\in A\big] \;d\red(x) =
\lambda\cdot \leb S\cdot  \P^* \big[(\red^*,\Psi^*)\in A\big].
\end{equation}
(Indeed this may be taken as a {\em definition} of the Palm
process). More generally, for any non-negative measurable $f$ on the
appropriate space,
\begin{equation}\label{palm2}
\E\int_{\R^d} f\big(\theta^{-x}(\red,\Psi),x\big) \;d\red(x) =
\lambda \int_{\R^d} \E^* f\big((\red^*,\Psi^*),x\big)\;dx.
\end{equation}

Let $\mat$ be a translation-invariant (one- or two- color) matching
scheme. If we let
\begin{equation}\label{new-x}
X:=|\mat^*(0)|
\end{equation}
denote the distance
from the origin to its partner under the Palm measure,
then \eqref{palm1} yields in particular
\begin{equation}\label{meaning-of-x}
\E \#\big\{x\in\rred\cap S: |x-\mat(x)|\leq r\big\}= \lambda \cdot
\leb S\cdot \P^*(X\leq r).
\end{equation}
Hence the above definition of $X$ is consistent
with the earlier notation in \eqref{def-x}.

Note that the tail bound \eqref{triv} for the Poisson process is
now an elementary consequence of \eqref{poisson-palm} and \eqref{new-x}.

\subsubsection*{Partial matching and mass transport}

A partial matching of a set $U$ is the edge set $m$ of a simple graph
$(U,m)$ in which each vertex has degree at most 1.  As before we write
$m(x)=y$ if $\{x,y\}\in m$, and in addition we write $m(x)=\infty$ if
$x$ is unmatched (i.e.\ has degree 0).  A {\dof one-color}
(respectively {\dof two-color}) {\dof partial matching scheme} $\mat$
is a point process on pairs which yields almost surely a partial
matching of $\rred$ (respectively between $\rred$ and $\bblue$).

\begin{prop}[Fairness]
\label{fair} Let $\red,\blue$ be simple point processes of finite
intensity, and let $\mat$ be a translation-invariant two-color
partial matching scheme of $\red$ and $\blue$. Then the process of
matched red points and the process of matched blue points have equal
intensity.
\end{prop}
\sloppy In particular, Proposition \ref{fair} shows that
translation-invariant perfect mat\-ching schemes are possible only
between two point processes of equal intensity. In addition, by
applying the result to the matching obtained by deleting all edges
longer than $r$, we see that $X$ is equal in law to the analogous
random variable defined in terms of a typical {\em blue}
point.\fussy

We prove Proposition \ref{fair} via the following lemma which will
be useful elsewhere.  (See \cite{blps} for background).

\begin{lemma}[Mass transport principle]\
\label{mt}
\begin{mylist}
\item Suppose $t:\Z^d\times\Z^d\to[0,\infty]$ satisfies $t(u+w,v+w)=t(u,v)$
for all $u,v,w\in\Z^d$, and write $t(A,B):=\sum_{u\in A,v\in
B}t(u,v)$.  Then
$$t(0,\Z^d)=t(\Z^d,0).$$
\item Suppose $T$ is a
random non-negative measure on $\R^d\times \R^d$ such that
$T(A,B):=T(A\times B)$ and $T(A+w,B+w)$ are equal in law for all
$w\in\Z^d$.  Then
$$\E T(Q,\R^d)=\E T(\R^d,Q).$$
\end{mylist}
\end{lemma}

\begin{proof}\

(i):
$t(0,\Z^d)=\sum_{u\in\Z^d}t(0,u)=\sum_{u\in\Z^d}t(-u,0)=t(\Z^d,0).$

(ii): Apply (i) to $t(u,v):=\E T(Q_u,Q_v)$.
\end{proof}

We sometimes call $t$ or $T$ a mass transport, and think of $t(A,B)$
or $T(A,B)$ as the amount of mass sent from $A$ to $B$.

\begin{proof}[Proof of Proposition \ref{fair}]
Apply Lemma \ref{mt} to the mass transport
$$T(A,B):=\#\{x\in\rred\cap A: \mat(x)\in B\}$$
in which each matched red point sends unit mass to its partner.
 Then $\E T(Q,\R^d)=\E\#\{x\in\rred\cap
Q:\mat(x)\neq\infty\}$, which is the intensity of matched red
points, while similarly $\E T(\R^d,Q)$ is the intensity of matched
blue points.
\end{proof}

\subsubsection*{Stable matching}

Following Gale and Shapley \cite{gale-shapley}, we say that a
partial matching $m$ of a set $U\subset\R^d$ is {\dof stable} if
there do {\em not} exist distinct points $x,y\in U$ satisfying
\begin{equation}\label{unstable}
|x-y|<\min\big\{|x-m(x)|,|y-m(y)|\big\},
\end{equation}
where $|x-m(x)|:=\infty$ if $x$ is unmatched.  A pair $x,y$ satisfying
\eqref{unstable} is called unstable.  (The motivation for this
definition is that each point prefers to be matched with closer
points, so an unstable pair $x,y$ prefer to divorce their current
partners and marry each other.)  Similarly, a partial bipartite
matching between two sets $U,V$ is called {\dof stable} if there do
not exist $x\in U$ and $y\in V$ satisfying \eqref{unstable}.

We call a set $U\subset \R^d$ {\dof non-equidistant} if there do
not exist $w,x,y,z\in U$ with $\{w,x\}\neq\{y,z\}$ and
$|w-x|=|y-z|>0$.  A {\dof descending chain} is an infinite
sequence $x_0,x_1,\ldots\in U$ for which the distances
$|x_i-x_{i+1}|$ form a strictly decreasing sequence.

\begin{prop}[Unique stable matching]\label{stable}
Let $\red$ be a translation-invari\-ant homogeneous point
processes of finite intensity.  (Respectively, let $\red,\blue$ be
point processes of equal finite intensity, jointly ergodic under
translations).  Suppose that almost surely $\rred$ (respectively
$\rred\cup\bblue$) is non-equidistant, and has no descending
chains. Then almost surely there is a unique stable partial
matching of $\rred$ (respectively between $\rred$ and $\bblue$).
Furthermore, it is almost surely a perfect matching, and it is
produced by the iterated mutually closest matching algorithm
described in the introduction.
\end{prop}

Under the conditions in Proposition \ref{stable}, the stable
matching has the following additional interpretation. Grow a ball
centered at each red point (respectively, each red point and each
blue point) simultaneously, so that at time $t$ all the balls have
radius $t$. Whenever two balls touch (respectively, whenever an
$\red$-ball and a $\blue$-ball touch), match their centers to each
other, and remove the two balls.

The conditions on the point processes in Proposition \ref{stable}
hold in particular for homogeneous Poisson processes, as proved in
\cite{haggstrom-meester}.  Clearly, under the conditions of the
proposition, the unique stable matching gives an
isometry-invariant factor matching scheme. We postpone the proof
of Proposition \ref{stable} to Section \ref{sec-stable}.

\section{Two-color matching}
%%%%%%%%%%%%%%%%%%%%%%%%%%%%

In this section we give proofs of Theorems \ref{2col-ub} and
\ref{2col-lb}.

\begin{proof}[Proof of Theorem \ref{2col-ub}(i)]
Use the stable matching (see Proposition \ref{stable} and Theorem
\ref{stab-2col}(ii)).
\end{proof}

\begin{proof}[Proof of Theorem \ref{2col-ub}(ii)]
We shall give a construction which works in all dimensions, and gives
a matching scheme with tails $\P^*(X>r)\leq C r^{-d/2}$.

First note that it is sufficient to give a translation-invariant
matching satisfying the required bound, for then we may obtain a
(randomized) isometry-invariant version by applying a uniformly
random isometry preserving the origin (i.e., chosen according to the
Haar measure on the compact group of such isometries) to
$(\red,\blue,\mat)$. Indeed, it suffices to give a matching scheme
which is invariant under translations by elements of $\Z^d$ and
which satisfies
\begin{equation}\label{tedious}
\E\#\{x\in\rred\cap Q:|x-\mat(x)|>r\}\leq C r^{-d/2},
\end{equation}
for then we may achieve a translation-invariant version satisfying
the same bound by similarly applying a uniformly random translation
in the unit cube $Q$.

We start by defining a sequence of successively coarser random
partitions of $\R^d$ into boxes in a $\Z^d$-invariant way.  Let
$\tau_0,\tau_1,\ldots$ be i.i.d.\ uniformly random elements of the
discrete cube $\{0,1\}^d$, independent of the point processes
$\red,\blue$. For each $k=0,1,\ldots$, define a {\dof $k$-box} to
be any subset of $\R^d$ of the form
$$[0,2^k)^d+2^k z+\sum_{i=0}^{k-1}2^i\tau_i,$$
where $z\in\Z^d$.

Now, given the point processes $\red,\blue$ and the partitioning
into boxes, define a matching as follows.  Within each $0$-box,
match as many red/blue pairs as possible in some arbitrary
pre-determined way.  (For definiteness, choose from among the
bipartite partial matchings of maximum cardinality the one which
minimizes the total edge length.) Remove those points which have
been matched. Now match as many red/blue pairs of the remaining
points as possible within each $1$-box, remove these matched
points, and repeat for $2$-boxes and so on. The union of all these
partial matchings clearly gives a $\Z^d$-invariant {\em partial}
matching scheme $\mat$ between $\red$ and $\blue$.

We shall prove that the partial matching $\mat$ satisfies
\eqref{tedious}. From this it follows by taking $r\to\infty$ that
almost surely every red point is matched, and hence by applying
Lemma \ref{fair} (to the $\R^d$-invariant version) that every blue
point is matched also.

Fix $k$, and call a red point {\dof $k$-bad} if it has not been
matched within its $k$-box by stage $k$ of the matching algorithm.
Suppose each $k$-bad red point distributes mass $1$ uniformly to
its $k$-box, i.e.\ let
$$T(A,B):=\sum_{x\in A\cap\rred:\text{$x$ is
$k$-bad}}2^{-dk}\leb\{y\in B: \text{$x$ and $y$ lie in the same
$k$-box}\}.$$ Then we have
\begin{multline}\label{mass-out}
\E T(Q,\R^d)=\E\#\{\text{$k$-bad red points in $Q$}\} \\
\geq
\E\#\{x\in \rred\cap Q: |x-\mat(x)|>2^k\surd d\}
\end{multline}
 (since a $k$-box has diameter $2^k\surd d$).

On the other hand, writing $W$ for the random $k$-box containing
$Q$,
\begin{multline}\label{mass-in}
\E T(\R^d,Q)=2^{-dk}\, \E \#\{\text{$k$-bad red points in $W$}\} \\
=2^{-dk}\, \E(\red(W)-\blue(W))^+= 2^{-dk}\, \E S^+,
\end{multline}
where $S:=\red[0,2^k)^d-\blue[0,2^k)^d$, since the location of $W$
is independent of $\red,\blue$.

The central limit theorem gives $\E S^+/\sqrt{2^{dk+1}}\to\E\chi^+$
as $k\to\infty$, where $\chi$ is a standard Gaussian.  Combining
this with \eqref{mass-out},\eqref{mass-in} and applying Lemma
\ref{mt} we deduce \eqref{tedious} for some $C=C(d)\in(0,\infty)$
and all $r=2^k\surd d$ with $k=0,1,2\ldots$. Hence by taking $2^k
\surd d \leq r < 2^{k+1} \surd d$ the same holds (with a modified
constant) for all $r>0$.
\end{proof}

\begin{proof}[Proof of Theorem \ref{2col-ub}(iii)]
We will deduce the result by a limiting argument from a result in
\cite{talagrand} on matchings of finite sets of points; a similar
argument was used in \cite{h-p-extra}.

The following is proved in \cite[equation (1.8)]{talagrand}.  Let
$d\geq 3$ and let $\red_n$ and $\blue_n$ each consist of $n$ point
masses whose locations are all independent and uniformly distributed
on $[0,1]^d$. Then for each $n$ there exists a two-color matching
scheme ${\cal F}_n$ between $\red_n$ and $\blue_n$ such that
\begin{equation}\label{tala}
 \P(G_n)
\geq 1-n^{-2}, \text{ where } G_n:=\bigg\{n^{-1} \sum_{x\in[\red_n]}
\exp \big(Cn|x-{\cal F}_n(x)|^d\big) \leq 2\bigg\}.
\end{equation}
 Here the constant $C$
depends on $d$ but not $n$.

Now let $\widetilde{{\cal F}_n}$ be ${\cal F}_n$ conditioned on the
event $G_n$ (this corrects a minor error in \cite{h-p-extra}).  We
construct a translation-invariant matching scheme ${\cal M}_n$ by
scaling $\widetilde{{\cal F}_n}$ to cover a cube of volume $n$, and
tiling $\R^d$ with identical copies of this matching, with the
origin chosen uniformly at random. More formally, regarding a
two-color matching scheme $\mat$ as a simple point process (i.e. a
random point measure) of {\em ordered} pairs $(x,y)\in
\R^d\times\R^d$ in which the presence of a point $(r,b)\in [\mat]$
indicates a matched pair $r\in\rred$ and $b\in \bblue$, we define
$$\mat_n(A\times B):=\sum_{z\in \Z^d}
\widetilde{{\cal F}_n}\Big(n^{1/d}(A+U+z)\times
n^{1/d}(B+U+z)\Big),$$ where $U$ is uniformly distributed on
$[0,1]^d$ and independent of $\widetilde{{\cal F}_n}$. Then
\eqref{tala} implies that for any Borel $A\subseteq\R^d$ we have
\begin{equation}\label{tight}
\E\iint \exp \big(C |x-y|^d\big) \ind_{x\in A} \mat_n(dx\times dy)
\leq 2 \leb A.
\end{equation}
(To check this, we first use invariance to deduce that the left side
must be a linear multiple of $\leb A$, and then take
$A=[0,n^{1/d}]^d$ to find the constant).

By \eqref{tight}, the random sequence $(\mat_n)$ is tight in the
vague topology of measures on $\R^d\times\R^d$ (see \cite[Lemma
16.15]{kallenberg}).  Therefore let $\mat$ be any subsequential
limit in distribution, and note that it has the following
properties.  It is a two-color matching scheme between the
marginal point processes $\mat(\cdot,\R^d)$ and
$\mat(\R^d,\cdot)$. These processes are independent Poisson point
processes of intensity 1 (this would clearly be true for any limit
constructed in the same way from the unconditioned matchings
${\cal F}_n$, therefore it holds for $\mat$ because $\P(G_n)\to
1$). The process $\mat$ inherits the translation-invariance of
$\mat_n$. Finally, it satisfies \eqref{tight} (with $\mat$
replacing $\mat_n$), which implies $\E^* e^{C X^d}\leq 2$ as
required.
\end{proof}

\paragraph{Remarks.}
The analogous results to Theorem \ref{2col-ub}(iii) for matchings of coin
flips on $\Z^d$ and for fair allocations of $\R^d$ may be proved by
following a similar limiting argument in the appropriate space (see
\cite{h-p-extra} for another variant).  The required matchings and
allocations on finite cubes exist by the results of \cite{talagrand}.

The result \eqref{tala} is a special case of a much more general
result in \cite{talagrand}, proved by deep (and indirect) methods.
A remark is made in \cite{talagrand} that the bound \eqref{tala}
can also be proved for an explicit matching obtained from the
construction in \cite{akt}.  Also, by results in \cite{grav}, in
$d\geq 3$ it is possible to obtain an explicit {\em ``transport''}
between two independent Poisson processes (i.e.\ a
translation-invariant random measure $T$ on $\R^d\times\R^d$ with
marginals $\red$ and $\blue$), satisfying a tail bound that is
exponential in a power of distance. (Specifically, construct the
``gravitational allocations'' -- see \cite{grav} -- for $\red$ and
$\blue$ independently, and let a red point $r$ send to a blue
point $b$ a mass equal to the volume of the intersection of the
cell of $r$ and the cell of $b$.)

We do not know the optimal tail behavior of translation-invariant
{\em factor} matching schemes of two independent Poisson processes
in dimensions $d\geq 2$.  We suspect that
there exist factor matchings satisfying the same bounds as
in the randomized case (see Theorem \ref{2col-ub}).  For the
analogous questions concerning coin flips on $\Z^d$, substantial
progress has been made by T. Soo \cite{soo} and A. Timar
\cite{timar-matching}.\parend

Theorem \ref{2col-lb} is Corollary 9 of \cite{h-p-extra}, which in
turn was deduced from results on ``extra head schemes'' in
\cite{h-liggett,liggett-tagged}.  Here we give a simple direct
argument.

\begin{proof}[Proof of Theorem \ref{2col-lb} (case d=1)]
Using \eqref{palm2} and Fubini's theorem, for $t>0$ we have
\begin{eqnarray*}
\lefteqn{\E\#\big\{x\in \rred\cap[0,2t]: \mat(x)\notin[0,2t]\big\}} \\
&&\leq  \E\#\big\{x\in \rred\cap[0,2t]: |\mat(x)-x|> x\wedge (2t-x)\big\} \\
&&= \int_0^{2t} \P^*\big[X>x\wedge(2t-x)\big] \;dx \\
%&= 2\int_0^t \P^*(X>x) \;dx \\
&&= 2\E^*(X\wedge t).
\end{eqnarray*}

Now the central limit theorem gives
$$\E\#\big\{x\in \rred\cap[0,2t]: \mat(x)\notin[0,2t]\big\}\geq
\E\big[(\red[0,2t]-\blue[0,2t])^+\big] \sim C t^{1/2}$$ as
$t\to\infty$ for some $C\in(0,\infty)$.  On the other hand, if $\E^*
X^{1/2}<\infty$ then the dominated convergence theorem gives
$t^{-1/2} \E^* (X\wedge t)\to 0$ as $t\to\infty$, so we obtain a
contradiction.
\end{proof}

For any $x,y\in\R^d$, we define the {\dof line-segment} $\seg{
x,y}:=\{\lambda x+ (1-\lambda)y:\lambda\in[0,1]\}\subset \R^d$.  The
following lemma will be used to derive a contradiction in the proof of
Theorem \ref{2col-lb} in the case $d=2$.

\begin{lemma}[Edge intersections]
\label{finite-intersection} For any translation-invariant 1-color or
2-color matching scheme $\mat$ (of any translation-invariant point
process(es)) which satisfies $\E^* X<\infty$, the number of matching
edges $\{x,y\}\in[\mat]$ such that the line segment $\seg{x,y}$
intersects the unit cube $Q$ has finite expectation.
\end{lemma}

\begin{proof}
Consider the mass transport
$$t(u,v):=\E\#\big\{x\in \rred\cap Q_u: \seg{ x,\mat(x)} \text{ intersects }
Q_v\big\}.$$ Since an edge of length $\ell$ intersects at most
$d(1+\ell)$ cubes of the form $z+Q$, where $z\in\Z^d$, we have
$$t(0,\Z^d)\leq d+d\,\E\sum_{x\in Q\cap\rred}
|x-\mat(x)|=d(1+\E^*X)<\infty,$$ hence by Lemma \ref{mt},
$$\E\#\{\text{matching edges intersecting }Q\}=t(\Z^d,0)<\infty.$$
\end{proof}

\begin{proof}[Proof of Theorem \ref{2col-lb} (case d=2)]
Without loss of generality we may assume that the matching scheme
$\mat$ is ergodic under translations of $\R^2$; if not we apply the
claimed result to the ergodic components.  Therefore suppose for a
contradiction that $\mat$ is an ergodic matching scheme satisfying
$\E^*X<\infty$.

For an ordered pair of distinct points $x,y\in\R^2$, we define the
random variable $K(x,y)$ to be the number of matching edges which
intersect the directed line segment $\seg{x,y}$ with the red point
on the left and the blue on the right. More precisely, $K(x,y)$ is
the number of pairs $\{r,b\}\in[\mat]$ with $r\in\rred$ and
$b\in\bblue$ such that $\seg{r,b}$ intersects $\seg{x,y}$, and
$\det\binom{ b-r}{y-x}>0$.

Lemma \ref{finite-intersection} and the assumption of finite mean
imply that $\E K(x,y)<\infty$ for any fixed $x,y$.  Note also the
additivity property that if $y\in\seg{x,z}$ (and these points
are deterministic), then a.s.\
$K(x,y)+K(y,z)=K(x,z)$. Fix a unit vector $u\in \R^2$.  Using
the ergodic theorem and the ergodicity of $\mat$ we deduce
$$\frac{K(0,nu)}{n} \xrightarrow{L^1} k(u) \quad\text{ as }n\to\infty,$$
where $k(u):=\E K(0,u)<\infty$. Since $K(v,v+nu)$ has the same law
as $K(0,nu)$, it follows that
\begin{equation}\label{ergodic}
\frac{K(v_n,v_n +nu)}{n} \xrightarrow{L^1} k(u) \quad\text{ as
}n\to\infty,
\end{equation}
for any deterministic sequence $v_n\in\R^2$.

Now denote the square $S:=[0,n]^2\subset\R^2$, and its corners
$s_1=(0,0)$; $s_2=(n,0)$; $s_3=(n,n)$; $s_4=(0,n)$.  The difference
$\red(S)-\blue(S)$ equals the number of edges crossing the boundary
of $S$ with the red point inside and the blue point outside minus
the number crossing in the opposite orientation. Hence,
\begin{equation}\label{compare}
\red(S)-\blue(S) = K_{+} - K_{-},
\end{equation}
where
$$K_+:=K(s_1,s_2)+K(s_2,s_3)+K(s_3,s_4)+K(s_4,s_1),$$
$$K_-:=K(s_1,s_4)+K(s_4,s_3)+K(s_3,s_2)+K(s_2,s_1).$$
(Note in particular that matching edges which intersect two sides of
$S$ do not contribute to the right side of \eqref{compare}, since
the contributions to $K_+$ and $K_-$ cancel.)

Writing $h=(1,0)$ and $v=(0,1)$, \eqref{ergodic} yields as
$n\to\infty$
$$\frac{K_+}{n}\xrightarrow{L^1} k(h)+k(v)+k(-h)+k(-v),$$
and
$$\frac{K_-}{n}\xrightarrow{L^1} k(v)+k(h)+k(-v)+k(-h),$$
hence $(K_+-K_-)/n \xrightarrow{L^1} 0$; that is $\E|K_+-K_-|=o(n)$.
On the other hand the central limit theorem gives
$\E|\red(S)-\blue(S)|=\Omega(n)$, contradicting \eqref{compare}.
\end{proof}

\section{One-color matching}
%%%%%%%%%%%%%%%%%%%%%%%%%%%%
\label{sec-1-col}

In this section we prove Theorems \ref{1col-1d} and \ref{trees}.

Consider the case $d=1$.  An {\dof adjacent} matching of a discrete
set $U\subset \R$ is one in which for every edge $\{x,y\}$, the
interval $(x,y)$ contains no points of $U$. Clearly for any given
infinite discrete $U$ there are exactly two adjacent matchings (one
in which the origin lies in $(x,y)$ for some edge $\{x,y\}$, and one
in which it does not).

\begin{proof}[Proof of Theorem \ref{1col-1d}(ii)]
Conditional on $\red$, choose one of the two adjacent matchings each
with probability $1/2$.  It is elementary to check that the
resulting matching scheme is isometry-invariant, and that
$\mat^*(0)$ is a symmetric two-sided exponential random variable.
\end{proof}

The following lemma states that, unsurprisingly, we cannot achieve
an adjacent matching without randomization.
  In what follows we say that
an edge $\{x,y\}$ of $\mat$ {\dof crosses} a site $z\in\R$ if
$z\in(x,y)$.
\begin{lemma}
\label{no-alt} Let $\red$ be a homogeneous Poisson process on $\R$.
There does not exist a one-color factor matching scheme where the
matching is almost surely adjacent.
\end{lemma}

\begin{proof}
Suppose on the contrary that $\mat$ is such a matching scheme.
Write ${\cal F}_S$ for the $\sigma$-algebra generated by the
restriction of $\red$ to $S\subseteq\R$.  Since the event
$$A:=\{0 \text{ is crossed by some edge}\}$$
lies in ${\cal F}_\R$, for every $\epsilon>0$ there exists
$r=r(\epsilon)<\infty$ and an event $A_\epsilon\in {\cal
F}_{[-r,r]}$ such that $\P(A\sym A_\epsilon)<\epsilon$.  Moreover,
by translation-invariance we can find $B_\epsilon\in {\cal
F}_{[-2r,0]}$ such that $\P(\{-r \text{ is crossed}\}\sym
B_\epsilon)<\epsilon$. For an adjacent matching, observing the
$\red$-points in an interval together with whether some
deterministic point in the
interval is crossed determines the matching on the interval a.s.
Hence
there exists $L_\epsilon\in {\cal F}_{[-2r,0]} \subset {\cal
F}_{(-\infty,0]}$ with $\P(A\sym L_\epsilon)<\epsilon$.
  Since this
is true for every $\epsilon>0$ we deduce that $A\in \overline{{\cal
F}_{(-\infty,0]}}$, where the bar denotes completion under $\P$.
% (for a proof, note that $\ind[L_{2^{-k}}]\to\ind[A]$ a.s.).
Similarly we have $A\in
\overline{{\cal F}_{[0,\infty)}}$, so $A$ is independent of itself
and has probability $0$ or $1$. But now it is easy to see that
neither of the two resulting matching schemes is
translation-invariant.
\end{proof}

\paragraph{Remark.}
Lemma \ref{no-alt} may be strengthened to the (equally
unsurprising) fact that the {\em only} a.s.\ adjacent matching
scheme is the one in the proof of Theorem \ref{1col-1d}(ii).
Indeed, for such a scheme consider the a.e.-defined function of
point configurations $f(\pi):=\P(O \text{ is crossed}\mid
\red=\pi)$ -- we must show that it equals $1/2$ a.e. By
translation-invariance, for a translation $\theta^{-t}$ we have
a.e.\ $f(\theta^{-t} \pi)\in\{f(\pi),1-f(\pi)\}$, according to the
parity of $\pi[0,t]$. Hence $|f(\pi)-1/2|$ is a
translation-invariant function of $\pi$, so by ergodicity of
$\red$ it is a.e.\ constant.  If the constant is $c\neq 1/2$ then
the translation-equivariant function $h_\pi(t):=f(\theta^{-t}\pi)$
assigns (a.e.\ with respect to $\pi$ and $t$) values $c$ and $1-c$
to alternate intervals between points of the configuration $\pi$,
and we can clearly use this to construct an adjacent factor
matching, contradicting Lemma \ref{no-alt}.
\parend

\begin{proof}[Proof of Theorem \ref{1col-1d}(i)]
Suppose that $\red$ is a Poisson point process of intensity $1$
and $\mat$ is a translation-invariant factor matching scheme. We
claim that
\begin{equation}\label{crossed}
\P(0 \text{ is crossed by infinitely many edges})=1.
\end{equation}

Suppose \eqref{crossed} is false. On the event that $0$ is crossed
by finitely many edges, a.s.\ the same is true for any other $r\in
\R$, because the difference between the number of edges crossing
$r$ and the number of edges crossing $0$ is at most the number of
red points between $0$ and $r$. By ergodicity, it then follows
that a.s.\ every $r\in \R$ is crossed only finitely many times. We
can now define a new matching scheme, by matching two adjacent red
points $x$ and $y$ if and only if the points $r$ between them are
crossed an odd number of times in the matching $\mat$. The new
matching is an adjacent factor matching, which contradicts
Lemma~\ref{no-alt} and proves~\eqref{crossed}.

To conclude, using \eqref{palm2} and Fubini's theorem we have
\begin{align*}
\E\#\{\text{edges crossing 0}\}
 &\leq\tfrac 12\E\#\big\{x\in\rred:|\mat(x)-x|>|x|\big\} \\
 &=\tfrac 12 \int_{-\infty}^{\infty} \P^*(X>|r|) \;dr \\
 &=\E^* X,
\end{align*}
so \eqref{crossed} implies $\E^* X=\infty$.
\end{proof}

\paragraph{Matching from a forest.}
The following construction will be used in the proof of Theorem
\ref{trees}.  Let $U$ be a countable infinite set, and suppose we
are given a locally finite forest with vertex set $U$ and one end
per tree (that is, a simple acyclic graph with finite degrees
where from each vertex there is exactly one singly-infinite
self-avoiding path). The {\dof parent} of vertex $x$ is the vertex
$y$ such that edge $(x,y)$ lies on the infinite path from $x$.  If
$y$ is the parent of $x$ then $x$ is a {\dof child} of $y$.
Suppose we are also given, for each vertex, a total ordering of
its children.

Under the above conditions, the following construction gives a
perfect matching of $U$. (Similar constructions were used in
\cite{ferrari-landim-thorisson} and \cite{h-p-trees}).  A {\dof leaf} is
a vertex with no children, and a {\dof twig} is a vertex which is
not a leaf but whose only children are leaves. Suppose $x$ is a
twig, and let $x_1,x_2,\ldots,x_k$ be its ordered children. Match
the pairs $\{x_1,x_2\},\{x_3,x_4\},\ldots$. If $k$ is odd, also
match $x_k$ to $x$.  Do this for all twigs.  Now remove those
vertices which have been matched, together with their incident
edges, and repeat the construction on the remaining graph. Repeat
indefinitely.

The above construction clearly gives a partial matching of $U$. To
see that it is a perfect matching, observe that for any given
vertex, at least one of its finitely many descendants is matched and
removed at every stage, so it is eventually matched itself.  (Vertex
$x$ is called a {\dof descendant} of $y$ if $y$ lies on the infinite
path from $x$.)

\sloppy Note the additional property that any matched pair were at
graph-theore\-tic distance at most 2 in the original forest. This
will enable us to derive tail bounds for the matching from tail
bounds on the forest. \fussy
\parend

Let $\red$ be a Poisson process of intensity 1.  The {\dof minimal
spanning forest} ${\cal F}$ of $\rred$ is the graph obtained from the
complete graph on $\rred$ by deleting every edge which is the longest
(in Euclidean distance) in some cycle.  In dimension $d=2$ it is known
\cite[p.~94]{alexander} that ${\cal F}$ is almost surely a locally
finite one-ended tree.  Now, ordering the children of each vertex $x$
by Euclidean distance from $x$ (say), and applying the above
construction, we get an isometry-invariant one-color factor matching
scheme.  This will enable us to prove Theorem \ref{trees} for $d=2$.
To prove the required tail bound we need the lemmas below.

Let $S:=\{y\in\R^d: |y|=1\}$ be the unit sphere. A {\dof cap} is a
proper subset of $S$ of the form $H=S\cap B(y,r)$, where $y\in S$. A
{\dof cone} of {\dof width} $w$ is a set of the form
$$V=V_H:=\big\{\alpha h: h\in H \text{ and } \alpha\in(0,\infty)\big\},$$
where $H$ is some cap of diameter $w$.

\begin{lemma}[Cones]
\label{cones} If $V=V_H$ is a cone of width $1$, then for any $x,y\in
V$,
$$|x|\geq |y| \text{ implies } |x-y|\leq |x|.$$
\end{lemma}

\begin{proof}
Let
$\beta:=|y|/|x|$, so that $|\beta x|=|y|$, and hence
$\beta x,y\in |y| H$.  Then
$$|x-y|\leq |x-\beta x|+|\beta x-y|
\leq (1-\beta)|x| + \diam(|y| H) = |x| -|y| + |y|.$$
\end{proof}

\begin{lemma}\label{msf}
In the minimal spanning forest ${\cal F}$ of a Poisson process
$\red$ of intensity 1 on $\R^d$,
$$\P^*(0 \text{ has an incident edge of length }>r)\leq C' e^{-C
r^d} \quad \forall r>0,$$ for $C,C'\in(0,\infty)$ depending only on
$d$.
\end{lemma}

\begin{proof}
For distinct $x,y\in\R^d$, define a set $S_{x,y}:=(V+x)\cap
B(x,|x-y|)$, where $V$ is some cone of width 1 such that $y\in V+x$.
For red points $x$ and $y$, if there exists another red point $z\in
S_{x,y}$ then $\cal F$ cannot have an edge from $x$ to $y$. This is
because $|x-z|<|x-y|$, while $|y-z|\le |y-x|$ by Lemma \ref{cones}, so
$\{x,y\}$ is the longest edge in the cycle $(x,y,z)$. Hence the
probability in question is at most
$$\P^*\Big[\red^*\big(V\cap B(0,r)\big)=0 \
\text{ for some cone $V$ of width 1}\Big].$$
Now let $H_1,\ldots,H_k$ be a finite subcover of the cover of the
unit sphere by caps of diameter $1/3$, and let $V_i:=V_{H_i}$ be the
associated cones of width $1/3$.  Any cone of width 1 must
completely contain at least one of the $V_i$, so the above
probability is at most
%%\begin{align*}
%&
\begin{multline*}
\P^*\Big[\red^*\big(V_i\cap B(0,r)\big)=0 \text{ for some
$i\in\{1,\ldots ,k\}$}\Big] \\
\leq %\; &
k \,\P^*\Big[\red\big(V_1\cap B(0,r)\big)=0\Big] \leq k
e^{-C r^d} .
\end{multline*}
%%\end{align*}
\end{proof}

\begin{lemma}\label{exponentials}
If non-negative random variables $Y,Z$ satisfy $\E
\,e^{cY^d}<\infty$ and $\E \,e^{cZ^d}<\infty$, then $\E
\,e^{c'(Y+Z)^d}<\infty$ where $c':=2^{-d}c$.
\end{lemma}

\begin{proof}
$\E \,e^{c'(Y+Z)^d}\leq \E \,e^{c(Y\vee Z)^d} = \E [e^{cY^d}\vee
e^{cZ^d}] \leq \E \,e^{cY^d}+\E \,e^{cZ^d}$.
\end{proof}

\begin{proof}[Proof of Theorem \ref{trees} (case $d=2$)]
Construct a matching $\mat$ from the minimum spanning forest ${\cal
F}$ as described above.  To prove the tail bound note that any
matched pair are either siblings or a parent and child. Hence if we
write $D(x)$ for the distance from $x$ to its parent in ${\cal F}$
then $|x-\mat(x)|\leq D(x)+D(\mat(x))$. Also, Lemma \ref{msf}
implies that $\E^* e^{C D(0)^2}<\infty$ for a certain constant
$C>0$.  Now defining
$$T(A,B):=\sum_{x\in \rred\cap A:\; \mat(x)\in B} e^{C D(x)^2},$$
the mass transport principle (Lemma \ref{mt}) gives $\E^* e^{C
D(\mat^*(0))^2}=\E^* e^{C D(0)^2}$, and by Lemma \ref{exponentials}
we obtain $\E^* e^{C' X^2}\leq \E^* e^{C'[D(0)+D(\mat^*(0))]^2}<\infty$.
\end{proof}

\paragraph{Remark.}  In dimensions $d\geq 3$ the minimal spanning
forest is believed to have one end per tree, but this is not
proved. Therefore we use a different forest (see below), which is
translation-invariant but not isometry-invariant.  Given
sufficient effort, it seems probable that a suitable one-ended
isometry-invariant forest could be constructed, giving an
isometry-invariant factor matching satisfying a similar bound.
\parend

\begin{proof}[Proof of Theorem \ref{trees} (case $d\ge 3$)]
Let $\red$ be a Poisson process of intensity 1.  As in the proof of
Theorem \ref{2col-ub}(ii), it suffices to give a
translation-invariant factor matching scheme, for then we may obtain
an isometry-invariant randomized version by applying a random
isometry preserving the origin.

The following construction is inspired by
\cite{ferrari-landim-thorisson}.  For $z=(z_1,z_2,\ldots ,z_d)\in
\R^d$ we write $\underline{z}:=(z_2,\ldots,z_d)\in\R^{d-1}$.  Define
the cone $V:=\{z\in \R^d: z_1>|\underline{z}|\}$.  Now for each red
point $x$, let $S(x)$ be the a.s.\ unique red point in $x+V$ for
which the first coordinate $S(x)_1$ is least, and put a directed
edge from $x$ to $S(x)$. Let $\mathcal G$ denote the resulting graph.

Since a.s.\ the out-degree of each vertex is 1 and there are no
oriented cycles, the graph ${\cal G}$ is clearly a forest.  The
mass-transport principle (Lemma \ref{mt}) shows that, under the
Palm measure, the expected in-degree of the origin is also $1$,
hence ${\cal G}$ is locally finite.  Furthermore, it is immediate
that
\begin{equation} \label{forest-bound}
\E^* e^{-C |S(0)|^d}<\infty
\end{equation}
 for some $C=C(d)>0$. We claim that
${\cal G}$ has one end per tree. Once this is proved, we can use it to
construct a matching as described above (ordering children by
distance), and the required tail bound may then be deduced from
\eqref{forest-bound} in the same way as in the above proof for the case $d=2$.

Turning to the proof that ${\cal G}$ has one end per tree, it suffices to
prove that ${\cal G}$ a.s.\ has no backward infinite path (that is,
no sequence of vertices and directed edges $\cdots \to x_2\to x_1
\to x_0$), for clearly from each vertex there is exactly one forward
infinite path.  Call a red point {\dof bad} if it lies on some
backward infinite path.  It is proved in \cite{alexander} that no
translation-invariant random forest in $\R^d$ can have a component with
more than 2 ends.  Assuming that bad points exist we shall obtain a
contradiction to this result.

Consider the `hyperplane' of cubes $L:=\{Q_u: u\in\Z^d, \,
u_1=-1\}$. If bad points exist then by invariance and ergodicity
there exist two fixed cubes $Q_u,Q_v\in L$ at distance at least $d+1$
from each other such that the event
$$A:=\{Q_u \text{ and } Q_v \text{ each contain a bad point}\}$$
has positive probability.  Note that whenever ${\cal G}$ has a
directed path from $x$ to $y$ then $y\in x+V$, or equivalently $x\in
y+(-V)$.  Hence for any red points $x\in Q_u$ and $y\in Q_v$, there
cannot be a directed path from $x$ to $y$ or from $y$ to $x$.
Recalling also that all out-degrees equal 1 we deduce that on $A$
there exist two {\em disjoint} backward infinite paths to red points
$x\in Q_u$ and $y\in Q_v$.  Also note that the event $A$ is
measurable with respect to the restriction of $\red$ to the
half-space ${\mathbb H}_-:=\R_-\times \R^{d-1}$.

Now we may construct an event $B$ of positive probability,
measurable with respect to the restriction of $\red$ to the
half-space ${\mathbb H}_+:=\R_+\times \R^{d-1}$, such that on $B$,
the forward infinite paths from any red points $x\in Q_u$ and
$y\in Q_v$ must coalesce.  Specifically, this will hold if a
sufficiently large region of ${\mathbb H}_+$ is empty except for
one red point which lies in $\bigcap_{z\in [(Q_u\cup Q_v)+V]\cap
{\mathbb H}_-}(z+V)$ (see Figure \ref{construct}). Now $A$ and $B$
are independent so $\P(A\cap B)>0$, but on the latter event,
${\cal G}$ has a component with at least 3 ends (formed by the
backward paths from the bad points in $Q_u$ and $Q_v$ together
with their joint forward path). This contradicts the result from
\cite{alexander} noted above.
\begin{figure}
\centering \resizebox{!}{3.2in}{\includegraphics{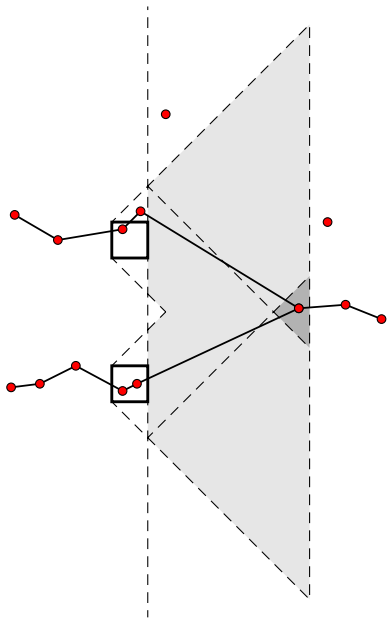}}
\caption{An event $B$ which forces coalescence of paths starting in
two given cubes: the light gray region should be empty while the
dark gray region should contain exactly one red point.} \label{construct}
\end{figure}
\end{proof}

\section{Stable Matching}
%%%%%%%%%%%%%%%%%%%%%%%%%
\label{sec-stable}

In this section we prove Theorems \ref{stab-1col} and
\ref{stab-2col} as well as Proposition \ref{stable}.  We start with
some relatively straightforward cases.

\begin{proof}[Proof of Theorem \ref{stab-1col}(ii)]
Call a red point $x$ {\dof $t$-bad} if $|\mat(x)-x|>t$, and note
that no two $t$-bad points may be within distance $t$ of each other,
for they would form an unstable pair.  Hence, using \eqref{palm1},
$$1\geq \E\#\big\{\text{$t$-bad red points in }B(0,\tfrac t2)\big\}
= \leb B(0,\tfrac t2) \cdot \P^*(X>t).$$
\end{proof}

\begin{proof}[Proof of Theorem \ref{stab-2col}(ii) (case $d=1$)]
Fix $t>0$.  We say that a red point $x$ is {\dof $t$-bad} if
$|\mat(x)-x|>t$.  Let $W$ be the set of $t$-bad red points in the
interval $[0,t]$, and denote the random variables
$$a=\min W;\quad b=\max W.$$
We claim that, provided $W\neq\emptyset$, every blue point in
$[a,b]$ is matched to a red point in $[a,b]$.  To prove this,
suppose on the contrary that $x\in[a,b]$ is a blue point matched
outside $[a,b]$. Without loss of generality suppose $\mat(x)<a$.
Then since $|\mat(x)-x|>|a-x|$ and $|\mat(a)-a|>t>|a-x|$, the pair
$a,x$ would be unstable, a contradiction.

Since elements of $W$ are $t$-bad, they cannot be matched within
$[a,b]\subseteq[0,t]$, so the above claim implies that
$\red[a,b]-\blue[a,b]\geq \#W$.  Hence
\begin{align*}
\#W&\leq \max_{[\alpha,\beta]\subseteq[0,t]}
\Big(\red[\alpha,\beta]-\blue[\alpha,\beta]\Big) \\
&\leq \max_{z\in[0,t]} F(z) - \min_{z\in[0,t]} F(z),
\end{align*}
where $F(z):=\red[0,z]-\blue[0,z]$. But $F$ is a continuous-time
simple symmetric random walk, so (by e.g.\ \cite[Proposition 13.13
and Theorem 14.6]{kallenberg}) taking the expectation of the above
inequality and using \eqref{palm1}, we obtain
$$t \,\P^*(X>t)=\E\#W\leq C \surd t$$
for a fixed constant $C\in(0,\infty)$.
\end{proof}

\begin{proof}[Proof of Theorem \ref{stab-2col}(ii) (case $d\geq 2$)]
See Theorem \ref{t.upperbd} at the end of this section.
\end{proof}

To prove the lower bounds Theorems \ref{stab-1col}(i) and
\ref{stab-2col}(i) we need the following simple properties of stable
matchings and Poisson processes.  The proofs are given after the
statement of Lemma \ref{local-mod}.

\begin{lemma}[Stable partial matchings]
\label{stable-det} Let $U$ (respectively, $U,V$) be (disjoint)
subset(s) of $\R^d$, and suppose that $U$ (respectively $U\cup V$)
is discrete, non-equidistant, and has no descending chains.  Then
there is a unique stable partial matching of $U$ (respectively
between $U$ and $V$), and it is produced by the iterated mutually
closest matching algorithm described in the introduction. In the
one-color case, there is at most one unmatched point, while in the
two-color case, all unmatched points are of the same color.
\end{lemma}
It should be noted that stable marriage problems do not in general
have unique solutions; see \cite{gale-shapley}.  The key to
uniqueness in our setting is that preferences are based on
distance, and are therefore symmetric.  Note also that some
condition on $U$ (or $U\cup V$) is needed in order to guarantee
existence and uniqueness of the stable matching.  For example, in
the one-color case with $d=1$, the set
$U=\{0\}\cup\{3^{-n}:n\in\Z^+\}$ has no stable partial matching,
while $U=\{\sum_{k=1}^n k^{-1}: n\in \Z^+\}$ has more than one.

\begin{lemma}[Modifications for 1-color stable matching]
\label{stable-mod}
Let $U\subset\R^d$ be a discrete, non-equidistant set with no
descending chains, and let $m$ be its unique stable partial
matching.
\begin{mylist}
\item If $\{x,y\}\in m$ is a matched pair then
$m\setminus\{\{x,y\}\}$ is the unique stable partial matching of
$U\setminus\{x,y\}$.
\item If $z\in\R^d\setminus U$ is such that
$U\cup\{z\}$ is non-equidistant, and $|m(x)-x|<|z-x|$ for all
$x\in U$, then $m$ is the unique stable partial matching of
$U\cup\{z\}$ (in particular, $z$ is unmatched).
\end{mylist}
\end{lemma}

\begin{lemma}[Monotonicity for 2-color stable matching]\label{mono}
Let $U,V,\{v'\}$ $\subset\R^d$ be disjoint sets, and suppose that
$U\cup V\cup\{v'\}$ is discrete, non-equidistant, and has no
descending chains.  Let $m$ be the stable bipartite partial
matching between $U$ and $V$, and let $m'$ be the stable bipartite
partial matching between $U$ and $V\cup\{v'\}$.  Then
$$|u-m'(u)|\leq |u-m(u)|\quad\forall u\in U,$$
(where as usual $|x-m(x)|:=\infty$ if $x$ is unmatched).
\end{lemma}
Lemma \ref{mono} states that adding an extra blue point makes the
matching no worse for red points.  Such results are well-known for
finite stable marriage problems -- see e.g.\
\cite{gale-shapley,knuth-stable}.

\begin{lemma}[Modifications for Poisson process]
\label{local-mod} Let $\red$ be a homogeneous Poisson process.
\begin{mylist}
\item
Let $U$ be a uniform random point in a set $S$ with $\leb
S\in(0,\infty)$, independent of $\red$.  The law of the point
process $\red+\delta_U$ obtained by adding a point at $U$ is
absolutely continuous with respect to the law of $\red$.
\item
Let ${\cal F}$ be a simple point process whose support $[{\cal F}]$
is a.s.\ a random {\em finite} subset of $\rred$.  The law of the
process $\red-{\cal F}$ obtained by removing the points of ${\cal
F}$ is absolutely continuous with respect to the law of $\red$.
\end{mylist}
\end{lemma}

\begin{proof}[Proof of Lemma \ref{stable-det}]
Consider the iterated mutually closest matching algorithm.
Non-equidistance ensures that it is well-defined.  We first claim that
every pair matched by the algorithm are matched to each other in every
stable partial matching.  This is proved by induction on the stage of
the algorithm: supposing the claim holds for all pairs matched so far,
any remaining mutually closest pair cannot be matched to points
removed earlier (by the inductive hypothesis) and they cannot be
matched further away than each other (by stability).

Now consider the set $N$ of points left unmatched by the
algorithm. We must show that these points are unmatched in every
stable partial matching. This is clear if $\#N\le 1$. Therefore,
suppose that $\#N\ge 2$ and consider first the one-color case. Let
$x_0\in N$, and let $x_1$ be the closest point to $x_0$ in
$N\setminus\{x_0\}$, (which exists, by discreteness). Inductively,
let $x_{n+1}$ be the point in $N\setminus\{x_n\}$ closest to
$x_n$. Clearly, $|x_n-x_{n+1}|\ge |x_{n+1}-x_{n+2}|$ for all $n\in
\mathbb N$. Since there is no descending chain, it follows that
there is some first $m\in\mathbb N$ such that $x_{m+1}=x_j$ with
$j\in\{0,1,\dots,m-1\}$. Because of non-equidistance, we must have
$j=m-1$. But then $x_{m-1}$ and $x_m$ are mutually closest in $N$,
which implies that they would have been matched by the algorithm
right after all other points had been removed from
$B\bigl(x_{m-1},|x_{m-1}-x_m|\bigr)\cup
B\bigl(x_m,|x_{m-1}-x_m|\bigr)$ (which would happen after a finite
number of stages by discreteness). This contradicts
$x_{m-1},x_m\in N$ and shows that $\#N\leq 1$.

In the two-color case, let $N_U$ be the unmatched points in $U$
and let $N_V$ be the unmatched points in $V$. If $N_U=\emptyset$,
then clearly all points in $N_V$ are unmatched in every stable
matching. The case $N_V=\emptyset$ is similar. If both $N_U$ and
$N_V$ are nonempty, we choose $x_0\in N_U$ and inductively let
$x_{2n+1}$ be the point in $N_V$ closest to $x_{2n}$ and let
$x_{2n+2}$ be the point in $N_U$ closest to $x_{2n+1}$. The
argument is then completed as in the one-color case.
\end{proof}

\begin{proof}[Proof of Proposition \ref{stable}]
Apply Lemma \ref{stable-det} and consider the random process ${\cal
N}$ of unmatched points in the unique stable partial matching - we
must show that $[{\cal N}]$ is almost surely empty.

In the one-color case, the lemma implies that $\#[{\cal N}]\le 1$.
  But if ${\cal N}$ has exactly one point with positive
probability then (after conditioning on this event) its location
would be a translation-invariant $\R^d$-valued random variable,
which is impossible.

In the two-color case, the lemma implies that $[{\cal N}]$ must be
empty, or consist entirely of red points or entirely of blue points.
By ergodicity, one of these three possibilities must have
probability 1.  But Lemma \ref{fair} implies that the processes of
unmatched red points and unmatched blue points have equal intensity,
so the latter two possibilities are ruled out.
\end{proof}

\begin{proof}[Proof of Lemma \ref{stable-mod}]
By Lemma~\ref{stable-det}, % Proposition \ref{stable},
we need only check that the claimed
matching is stable.  In (i) this is immediate, since any unstable
pair would have been unstable in the original matching.  Similarly
in (ii), the given condition ensures that $z$ does not form an
unstable pair with any $x\in U$.
\end{proof}

\begin{proof}[Proof of Lemma \ref{mono}]
Suppose on the contrary that for some $u\in U$ we have
$|u-m'(u)|>|u-m(u)|$.  In particular $m(u)\neq\infty$, so write
$v:=m(u)\in V$.  Stability of $(u,v)$ in $m'$ implies $|v-m'(v)|\leq
|v-u|$, but $m'(v)\neq u$ so non-equidistance implies that the
previous inequality is strict; we write $u_1:=m'(v)\in U$.
Similarly, by stability of $(u_1,v)$ in $m$ we have
$|u_1-m(u_1)|<|u_1-v|$; write $v_1:=m(u_1)\in V$.  Iterating this
argument gives a descending chain $u,v,u_1,v_1,u_2,v_2,\ldots$.
\end{proof}

\begin{proof}[Proof of Lemma \ref{local-mod}]
(i):  It is elementary to check that the Radon-Nikodym derivative of
the laws is $p(N-1)/p(N)$, where $N:=\red(S)$, and $p(k)$ is the
probability that $N=k$.

(ii): Let $A$ be some measurable set such that $\P(\mathcal
R-\mathcal F\in A)>0$. We need to show that $\P(\mathcal R\in
A)>0$. Since a.s.\ $[\mathcal F]$ is finite and $[\mathcal R]$ is
discrete, there is a.s.\ a finite random set of balls with
rational centers and radii such that $[\mathcal F]$ is the
intersection of $[\mathcal R]$ with the union of these balls.
 Therefore, there is a deterministic finite union of open balls
$W$ such that $\delta:=\P(\mathcal R-\mathcal F\in A,\,[\mathcal
R]\cap W=[\mathcal F])>0$. Let $\mathcal R_1$ denote the
restriction of $\mathcal R$ to the complement of $W$, and let
$A_1$ be the event that $\P\bigl(\mathcal R-\mathcal F\in
A,\,[\mathcal R]\cap W= [\mathcal F]\bigm|\mathcal
R_1\bigr)>\delta/2$.  Note that $\P(A_1)\geq \delta/2$.  On the
event $A_1$, with positive probability we have $\mathcal
R-\mathcal F=\mathcal R_1$ and therefore $\mathcal R_1\in A$.  But
$A_1$ is $\sigma(\mathcal R_1)$-measurable, so we must have
$\mathcal R_1\in A$ a.s.\ on $A_1$.  Since $\mathcal R_1$ and
$\mathcal R(W)$ are independent, we deduce
$$
\P(\mathcal R\in A)\ge \P\big(A_1,\,\mathcal R(W)=0\big)=
\P(A_1)\,\P\big(\mathcal R(W)=0\big)>0.
$$
\end{proof}

We now turn to the proofs of the lower bounds.

\begin{proof}[Proof of Theorem \ref{stab-1col}(i)]
Let $\mat$ be the one-color stable matching, and consider the random
set
\begin{equation}\label{set-e}
H=H(\red):=\{x\in\rred: |x-\mat(x)|>|x|-1\}.
\end{equation}
 This is
the set of red points that would prefer some red point in the unit
ball $B(0,1)$ (if one were present in the correct location) over
their current partners.  We claim that
\begin{equation}\label{theclaim}
\P(\#H=\infty)=1.
\end{equation}
Once this is proved, we obtain the required result as follows,
using \eqref{palm2} and Fubini's theorem:
\begin{align*}
\infty=\E\#H&=\int_{\R^d} \P^*(X>|x|-1)\;dx \\
&= \int_0^\infty \P^*(X+1>t)\, c\, t^{d-1} \;dt \\
&= \frac{c}{d}\, \E^*\bigl[(X+1)^d\bigr],
\end{align*}
hence $\E^* X^d=\infty$.

Returning to the claim \eqref{theclaim}, suppose on the contrary
that $H$ is finite with positive probability.  For each point
configuration $\red$, construct a modified configuration
$\widehat{\red}$ as follows:
\begin{mylist}
 \item if $\# H<\infty$, remove all the points in $H\cup \{\mat(x):x\in H\}$;
 \item add a uniformly random point in $B(0,1)$, independently of
 $\red$.
\end{mylist}
Using Lemma \ref{local-mod}, the law of the random configuration
$\widehat{\red}$ is absolutely continuous with respect to that of
$\red$. Now by Lemma \ref{stable-mod}, whenever $\#
H(\red)<\infty$, the stable partial matching of $[\widehat{\red}]$
has an unmatched point (the one added in (ii)), hence this happens
with positive probability. Absolute continuity therefore implies
that with positive probability the stable partial matching of
$\rred$ has an unmatched point, contradicting Proposition
\ref{stable}.
\end{proof}

\begin{proof}[Proof of Theorem \ref{stab-2col}(i)]
Define the random set $H$ exactly as in \eqref{set-e} (now it is the
set of red points which would prefer a blue point in $B(0,1)$). We
will prove that $\P(\#H=\infty)=1$, whereupon the result follows as
in the proof of Theorem \ref{stab-1col}(i).

Fix any $k<\infty$; we will prove that $\P(\#H \geq k)=1$.  Let
$\blue'$ be obtained from $\blue$ by adding $k$ independent
uniformly random points in $B(0,1)$. By Lemma \ref{local-mod}(i),
the law of $(\red,\blue')$ is absolutely continuous with respect to
that of $(\red,\blue)$. Hence, by Proposition \ref{stable}, almost
surely all the $k$ added points are matched in the stable matching
between $\rred$ and $[\blue']$. By Lemma \ref{mono}, it follows that
the partners of the added points were matched as far away or further
in the stable matching with $\bblue$, so these partners lie in $H$,
and thus $\#H \geq k$ as required.
\end{proof}

\paragraph{Remark.} As stated in the introduction, Theorem~\ref{stab-2col}
holds also for the stable matching of heads (red) to tails (blue)
on $\Z^d$ (given some tie-breaking rule). In order to adapt the
proof of (i) to that setting, we claim that the set $H$ of red
sites $v\in\Z^d$ which would prefer the origin to their current
partner (if the origin was blue)
 must be infinite. Indeed, if $H$ is finite,
then a contradiction is obtained by considering the configuration in
which the sites in $H$ are recolored blue and the origin is
colored blue.
\medskip

Finally we prove the upper bound for the two-color stable matching
in $d\geq 2$.

\begin{thm}\label{t.upperbd}
In the two-color stable matching of two independent Poisson
processes of intensity 1 in $\R^d$, $d\ge 2$, we have
\begin{equation}
\label{e.tailbd} \P^*(X>r)\le C\,r^{-s}\,,
\end{equation}
where $C=C(d)\in(0,1)$, and $s=s(d)$ is the unique solution in
$(0,1)$ of the equation
\begin{equation}
\label{e.seq} 2\,\omega_d\int_1^2(t-1)^{d-1}\,t^{-s}\,dt =
\frac{\omega_{d-1}}{d-1}\int_0^2\bigl(1-(t/2)^2\bigr)^{\frac{d-1}2}\,t^{-s}\,dt\,,
\end{equation}
and $\omega_d$ denotes the $(d-1)$-dimensional volume of the unit
sphere in $\R^d$.
\end{thm}

For $d=2$,~\eqref{e.seq} simplifies to
$$
\sqrt{\pi}\,(2^s-2\,s)\,\Gamma\bigl((2-s)/2\bigr)=\Gamma\bigl((3-s)/2\bigr)\,,
$$
and for general $d$ the integrals can be evaluated in terms of
hypergeometric functions.  The numerical values (rounded to the
nearest $10^{-3}$) of $s$ at $d=2$, $3$, $10$ and $100$ are $0.496$,
 $0.449$, $0.339$ and $0.224$, respectively. It is not hard to see that
$s_d\,\log d$ stays bounded away from $0$ and $\infty$ as
$d\to\infty$.

\begin{proof}[Proof of Theorem \ref{t.upperbd}]
Set $\alpha(r):= \P^*(X>r)$. Fix some $R>0$ and consider the ball
$B=B(0,R)$ of radius $R$ about $0$.  Set
$$
Y_\red:=\bigl\{x\in\rred\cap B:|x-\mat(x)|> R+|x|\bigr\},
$$
and similarly
$$
Y_\blue:=\bigl\{x\in\bblue\cap B:|x-\mat(x)|> R+|x|\bigr\}.
$$
First, observe that if $x\in Y_\red$, then $\mat(x)\notin B$. Next,
note that if $x\in Y_\red$, then $x$ prefers any blue point in $B$
to its partner. Since the corresponding statement also holds for
$Y_\blue$, we have
\begin{equation}\label{e.principal}
\text{if $Y_\red\ne\emptyset$ then $Y_\blue=\emptyset$}.
\end{equation}
(This is the principal observation on which the proof rests.) Let
$Q_\blue$ denote the set of blue points in $B$ that are matched
outside of $B$, and similarly for $Q_\red$.

Let $Z:=\red (B)-\blue(B)$, and note that $\#Q_\red-\#Q_\blue = Z$.
Therefore, $\#Y_\red \le \#Q_\blue+Z$, and~\eqref{e.principal} gives
$$
\#Y_\red+\#Y_\blue \le \#Q_\blue+Z^+.
$$
Our bound will follow by taking the expectation of both sides of
this inequality. Since $\E\#Y_\red=\E\#Y_\blue$ and $\E(Z^+)\le
C\,R^{d/2}$ for some fixed constant $C=C_d$ (which may depend only
on $d$), we get
\begin{equation}\label{e.thefunbegins}
2\,\E\#Y_\red \le C\,R^{d/2}+\E\#Q_\blue\,.
\end{equation}
By \eqref{palm2}, it is easy to express the left hand side in terms
of $\alpha$, namely,
\begin{equation}\label{e.Ybd}
\begin{aligned}
\E\#Y_\red &= \int_B \alpha(R+|x|)\,dx
\\&
= \omega_{d} \int_0^R \alpha(R+r)\, r^{d-1} \,dr = \omega_{d}
\int_R^{2R} \alpha(r)\, (r-R)^{d-1} \,dr \,.
\end{aligned}
\end{equation}

The proof will proceed by expressing $\E\#Q_\blue$ in terms of
$\alpha$ and using \eqref{e.thefunbegins}. Before embarking on the
full argument we note the following simplified version which already
gives a power law upper bound on $\alpha$. If a blue point in $B$ is
matched outside $B$ then the length of its edge is at least its
distance to the boundary of $B$, hence \eqref{palm2} gives
\begin{equation}\label{e.easy}
\begin{aligned}
\E\#Q_\blue &\le \int_B \alpha(R-|x|)\,dx
\\&
= \omega_{d} \int_0^R \alpha(R-r)\, r^{d-1} \,dr = \omega_{d}
\int_0^{R} \alpha(r)\, (R-r)^{d-1} \,dr \,.
\end{aligned}
\end{equation}
Substituting \eqref{e.Ybd} and \eqref{e.easy} into
\eqref{e.thefunbegins} and using the fact that $\alpha$ is
decreasing yields a bound for $\alpha(2R)$ in terms of $\alpha(r)$
for $r\in[0,R]$, and it is straightforward to deduce (by induction
on $k$) that $\alpha(2^k)\leq C' (2^k)^{-s'}$ for some
$C'=C'(d)\in(0,\infty)$ and $s'=s'(d)\in(0,1)$.

In order to get a better power we will instead use an exact
expression for $\E\#Q_\blue$, and analyze the resulting inequality
more carefully.  Denote the unit sphere
$S^{d-1}:=\{z\in\R^d:|z|=1\}$.  The intensity of the process of
pairs $(x,u)\in\R^d\times S^{d-1}$ such that $x\in[\blue]$,
$|x-\mat(x)|>r$ and $(\mat(x)-x)/|\mat(x)-x|=u$ is precisely
$\alpha(r)/\omega_{d}$. (That is, the expected number of such pairs
in any set $A\subset \R^d\times S^{d-1}$ is $\alpha(r)/\omega_{d}$
times the volume of $A$.) For $x\in B$ and $u\in S^{d-1}$, let
$q(x,u):=\inf\{t\ge 0: x+t\,u\notin B\}$, and fix some $u_0\in
S^{d-1}$. Then
\begin{equation}
\label{e.Qb} \E \#Q_\blue = \frac 1{\omega_d} \int_{S^{d-1}}\int_B
\alpha\bigl(q(x,u)\bigr) \,dx\,du = \int_B
\alpha\bigl(q(x,u_0)\bigr) \,dx \,,
\end{equation}
where the last equality is a consequence of rotational symmetry.
Let $L$ denote the orthogonal projection onto the subspace of
$\R^d$ orthogonal to $u_0$; that is, $L\,z= z-(z\cdot u_0)\,u_0$.
For $x\in B$ define $f(x)=L\,x+q(x,u_0)\,u_0$. Since
$L\,f(x)=L\,x$ and $f(x+t\,u_0)=f(x)-t\,u_0$, it follows by
differentiation that $f$ is measure preserving. This allows us to
use the substitution $z=f(x)$ and write
$$
\E \#Q_\blue = \int_{f(B)} \alpha\bigl(z\cdot u_0\bigr) \,dz =
\int_0^{2R} \mu_{d-1}\bigl\{z\in f(B): z\cdot u_0
=r\bigr\}\,\alpha(r)\,dr \,,
$$
where $\mu_{d-1}$ is $(d-1)$-dimensional Lebesgue measure and the
last equality follows by Fubini. Now,
$$
\mu_{d-1}\bigl\{z\in f(B): z\cdot u_0 =r\bigr\}=
\mu_{d-1}\bigl\{L\,x: x\in B,\,q(x,u_0)=r\bigr\}\,.
$$
Note that the set $ \{L\,x: x\in B,\,q(x,u_0)=r\bigr\}$ is precisely
the set of sites $z\in L\,\R^d$ such that $z- (r/2)\,u_0\in B$,
which is $\bigl\{z\in L\,\R^d: |z|< \sqrt{R^2-(r/2)^2}\bigr\}$. The
$(d-1)$-volume of this set is just
$\bigl(R^2-(r/2)^2\bigr)^{(d-1)/2}$ times the volume of the
$(d-1)$-dimensional unit ball. Since the volume of the
$(d-1)$-dimensional unit ball is $\omega_{d-1}/(d-1)$, we get
\begin{equation}\label{e.outedges}
\E \#Q_\blue = \frac{\omega_{d-1}}{d-1}\,\int_0^{2R}
\bigl(R^2-(r/2)^2\bigr)^{(d-1)/2}\,\alpha(r)\,dr\,.
\end{equation}

Now, taking into account~\eqref{e.thefunbegins}, \eqref{e.Ybd} and
\eqref{e.outedges}, we obtain
$$
\int_0^{\infty} g(r/R)\,\alpha(r)\,dr \le C\, R^{1-d/2},
$$
where
$$
g(t):= 2\,\omega_d\,\ind_{[1,2]}(t)\, (t-1)^{d-1}-
\frac{\omega_{d-1}}{d-1}\, \ind_{[0,2]}(t)
\,\bigl(1-(t/2)^2\bigr)^{\frac{d-1}2}.
$$
This bound will be useful when $R$ is large. For $R$ small, we use
the trivial estimate
$$
\int_0^\infty g(r/R)\,\alpha(r)\,dr\le \int_0^\infty g(r/R)\,dr \le
\int_0^{2R} \|g\|_\infty\,dr = 2\,R\,\|g\|_\infty\,.
$$
Combining these two estimates, we get
\begin{equation}
\label{e.ga} \int_0^\infty g(r/R)\,\alpha(r)\,dr\le
\min\{C\,R^{1-d/2},2\,R\,\|g\|_\infty\}\,.
\end{equation}
We will get our desired bound on $\alpha$ by taking an appropriate
average of~\eqref{e.ga} with respect to $R$.

Note that the set of $s$ satisfying $\int_0^\infty
g(t)\,t^{-s}\,dt=0$ is precisely the set of $s$
satisfying~\eqref{e.seq}. Observe that $g(t)$ is supported on
$[0,2]$ and is continuous and monotone increasing there. Moreover,
$g(0)<0<g(2)$. Therefore, there is some $s=s_d\in (-\infty,1)$
such that $\int_0^\infty g(t)\,t^{-s}\,dt=0$. We claim that $s$ is
unique. Indeed, let $s'<s$ and let $t_0$ be the unique solution of
$g(t)=0$ in $(0,2)$. Then $t_0^{s'-s}\,t^{-s'}< t^{-s}$ precisely
when $t<t_0$. Therefore $t_0^{s'-s}\int_0^\infty
g(t)\,t^{-s'}\,dt> \int_0^\infty g(t)\,t^{-s}\,dt=0$, proving
uniqueness. Next, we claim that $s>0$. Observe that if we replace
$\alpha$ by $1$ we get the volume of $B$ (that is,
$R^d\,\omega_d/d$) in~\eqref{e.Ybd} and~\eqref{e.Qb}. The
algebraic manipulations within and following these equalities are
valid for any measurable bounded function in place of $\alpha$.
Therefore $\int_0^\infty g(t)\,dt = \mu_d(B(0,1))=\omega_d/d>0$,
which implies $s>0$.

Since $\int_0^\infty g(t)\,t^{-s}\,dt=0$, a change of variables
$t=r/\rho$ gives
\begin{equation}
\label{e.geq}
  \int_0^\infty g(r/\rho)\,\rho^{s-2}\,d\rho=0\,.
\end{equation}
Set
$$
G_R(r):=\int_0^R g(r/\rho)\,\rho^{s-2}\,d\rho\,.
$$
We claim that $G_R(r)\ge 0$ for all $r>0$, and that
$$
C_0:=\inf\bigl\{G_R(r)\,R^{1-s}: R>0,\,r\in[R/2,R]\bigr\}>0\,.
$$
As before, let $t_0$ be the unique solution of $g(t)=0$ in
$(0,2)$. If $r/R\ge t_0$, then $g(r/\rho)\ge 0$ for $\rho<R$, and
hence $G_R(r)\ge 0$. On the other hand, if $r/R<t_0$, then
$g(r/\rho)\le 0$ for all $\rho>R$ and~\eqref{e.geq} gives
$G_R(r)=-\int_R^\infty g(r/\rho)\,\rho^{s-2}\,d\rho\ge 0$. Since
$g(t)>0$ for $t\in(t_0,2)$ and $g(t)<0$ for $t\in(0,t_0)$, the
above reasoning actually gives $G_R(r)>0$ for $r\in(0,2\,R)$.
Since $G_R$ is continuous, this implies $\inf_{r\in[1/2,1]}
G_1(r)>0$. A change of variables gives $G_R(r)=
R^{s-1}\,G_1(r/R)$, which now proves $C_0>0$.

From the monotonicity of $\alpha$, the definition of $C_0$ and from $G_R(r)\ge 0$,
we now get
\begin{equation}
\label{e.alphabd} C_0\,R^s\,\alpha(R)/2 \le C_0\,R^{s-1}\int_{R/2}^R
\alpha(r)\,dr\le \int_0^\infty G_R(r)\,\alpha(r)\,dr \,.
\end{equation}
Note that
$$
\int_0^R\int_0^\infty \rho^{s-2} \,|g(r/\rho)|\,\alpha(r)\,dr\,d\rho
\le \int_0^R\int_0^\infty \rho^{s-2} \,\|g\|_\infty\,1_{r\le 2\rho}
\,dr\,d\rho <\infty\,.
$$
Therefore, Fubini and~\eqref{e.ga} give
\begin{equation}
\label{e.final}
\begin{aligned}
\int_0^\infty G_R(r)\,\alpha(r)\,dr & =\int_{0}^R\left(\rho^{s-2}\,
\int_{0}^\infty g(r/\rho)\,\alpha(r)\,dr\right)\,d\rho
\\&
\le \int_{0}^R\rho^{s-2}\,
\min\{C\,\rho^{1-d/2},2\,\rho\,\|g\|_\infty \} \,d\rho
\\&
\le \int_{1}^R C\,\rho^{s-1-d/2}\,d\rho + 2\,\|g\|_\infty \int_0^1
\rho^{s-1} \,d\rho \,.
\end{aligned}
\end{equation}
Since the right hand side is bounded in $R$ (this is where we use
$d>1$), this and~\eqref{e.alphabd} imply~\eqref{e.tailbd}.
\end{proof}

\paragraph{Remarks.}
In order to adapt the proof of Theorem \ref{t.upperbd} to the
stable allocation of Lebesgue to Poisson, we replace $\#Y_\red$
with the volume of sites $z\in B$ whose Poisson point is at
distance greater than $R+|z|$, and replace $\#Y_\blue$ with the
sum over Poisson points $x\in B$ of the volume of $x$'s territory
that is at distance greater than $R+|x|$. The mass transport
principle easily shows that these two quantities have the same
expectation. A similar remark applies to $\#Q_\red$ and
$\#Q_\blue$. This allows us to obtain the analog
of~\eqref{e.thefunbegins}.

To adapt the proof to the setting of a stable matching of heads to
tails in $\Z^d$, we apply a uniform random translation in
$[0,1)^d$, and then apply a random isometry preserving the origin.
Then the law of the matching is invariant under isometries of
$\R^d$, and the above proof applies.

\section*{Open Problems}

\begin{mylist}
\item For the two-color stable matching of two independent Poisson
processes, what is the correct power law for the tail behavior of
$X$ in dimensions $d\geq 2$?  We conjecture that $\E^* X^\alpha
<\infty$ if and only if $\alpha<d/2$.

\item Does there exist a translation-invariant matching of two
independent Poisson processes in $\R^2$ such that the line
segments connecting matched pairs do not cross?

\fussy
\end{mylist}

\bibliography{bib}

\section*{ }

Alexander E. Holroyd: {\tt holroyd(at)math.ubc.ca} \\
University of British Columbia, 121-1984 Mathematics Rd, \\
Vancouver BC V6T 1Z2, Canada

\vspace{4mm} \noindent Robin Pemantle: {\tt pemantle(at)math.upenn.edu} \\
David Rittenhouse Laboratories, 209 S 33rd St, \\
Philadelphia PA 19104, USA

\vspace{4mm} \noindent Yuval Peres: {\tt peres(at)microsoft.com} \\
Oded Schramm: \\
Microsoft Research, One Microsoft Way, \\
Redmond WA 98052, USA

\end{document}